\documentclass[11pt,reqno]{amsart}
\usepackage{import,hyperref}
\usepackage[latin1]{inputenc}

\newtheorem{theorem}{Theorem}[section]
\newtheorem{proposition}[theorem]{Proposition}
\newtheorem{conjecture}[theorem]{Conjecture}
\newtheorem{lemma}[theorem]{Lemma}
\newtheorem{corollary}[theorem]{Corollary}

\theoremstyle{definition}
\newtheorem{definition}[theorem]{Definition}

\usepackage{amssymb,amsmath,amsthm,amsfonts}
\usepackage{mathtools}

\usepackage{eucal,latexsym}
\usepackage{ytableau}
\usepackage{shuffle}
\usepackage{commath} % \abs

\newcommand{\Z}{\mathbb{Z}}

\newcommand{\ra}{\rightarrow}

\newcommand{\sseq}{\subseteq}

\newcommand{\fk}{\mathfrak}
\DeclareMathOperator{\perm}{perm}
\DeclareMathOperator{\wt}{wt}

\DeclareMathOperator{\ASM}{ASM}

\DeclareMathOperator{\BPD}{BPD}
\DeclareMathOperator{\bpd}{bpd}
\DeclareMathOperator{\mBPD}{mBPD}
\DeclareMathOperator{\mbpd}{mbpd}

\usepackage{tikz}

\newcommand{\smalltilescale}{.2}
\newcommand{\tilescale}{.27}
\newcommand{\bigtilescale}{.5}

\newcommand{\pipecolor}{UForange}

\newcommand{\smallblank}{\blankTEMPLATE{\smalltilescale}}

\newcommand{\smalljelbow}{\jelbowTEMPLATE{\smalltilescale}}

\newcommand{\blank}{\blankTEMPLATE{\tilescale}}
\newcommand{\horizontal}{\horizontalTEMPLATE{\tilescale}}
\newcommand{\vertical}{\verticalTEMPLATE{\tilescale}}
\newcommand{\cross}{\crossTEMPLATE{\tilescale}}
\newcommand{\bump}{\bumpTEMPLATE{\tilescale}}
\newcommand{\relbow}{\relbowTEMPLATE{\tilescale}}
\newcommand{\jelbow}{\jelbowTEMPLATE{\tilescale}}

\newcommand{\bigblank}{\blankTEMPLATE{\bigtilescale}}
\newcommand{\bighorizontal}{\horizontalTEMPLATE{\bigtilescale}}
\newcommand{\bigvertical}{\verticalTEMPLATE{\bigtilescale}}
\newcommand{\bigcross}{\crossTEMPLATE{\bigtilescale}}
\newcommand{\bigbump}{\bumpTEMPLATE{\bigtilescale}}
\newcommand{\bigrelbow}{\relbowTEMPLATE{\bigtilescale}}
\newcommand{\bigjelbow}{\jelbowTEMPLATE{\bigtilescale}}

\newcommand{\DRAWGRID}{
    \draw[step=1.0,black,opacity=0.5,thin,xshift=0.5cm,yshift=0.5cm] (0,0) grid (1,1);
}

\newcommand{\blankTEMPLATE}[1]{
    \begin{tikzpicture}[scale=#1]
    \DRAWGRID
    \end{tikzpicture}
}
\newcommand{\horizontalTEMPLATE}[1]{
    \begin{tikzpicture}[scale=#1]
    \DRAWGRID
    \draw[color=\pipecolor,very thick]
    (0.5,1) to (1.5,1);
    \end{tikzpicture}
}
\newcommand{\verticalTEMPLATE}[1]{
    \begin{tikzpicture}[scale=#1]
    \DRAWGRID
    \draw[color=\pipecolor,very thick]
    (1,0.5) to (1,1.5);
    \end{tikzpicture}
}
\newcommand{\crossTEMPLATE}[1]{
    \begin{tikzpicture}[scale=#1]
    \DRAWGRID
    \draw[color=\pipecolor,very thick]
    (1,0.5) to (1,1.5)
    (0.5,1) to (1.5,1);
    \end{tikzpicture}
}
\newcommand{\bumpTEMPLATE}[1]{
    \begin{tikzpicture}[scale=#1]
    \DRAWGRID
    \draw[color=\pipecolor,very thick]
    (1,0.5) to (1,0.75) to[bend left] (1.25,1) to (1.5,1)
    (0.5,1) to (0.75,1) to[bend right] (1,1.25) to (1,1.5);
    \end{tikzpicture}
}
\newcommand{\relbowTEMPLATE}[1]{
    \begin{tikzpicture}[scale=#1]
    \DRAWGRID
    \draw[color=\pipecolor,very thick]
    (1,0.5) to (1,0.75) to[bend left] (1.25,1) to (1.5,1);
    \end{tikzpicture}
}
\newcommand{\jelbowTEMPLATE}[1]{
    \begin{tikzpicture}[scale=#1]
    \DRAWGRID
    \draw[color=\pipecolor,very thick]
    (0.5,1) to (0.75,1) to[bend right] (1,1.25) to (1,1.5);
    \end{tikzpicture}
}

\definecolor{UForange}{RGB}{250,70,22} % #FA4616
\definecolor{UFblue}{RGB}{0,33,165} % #0021A5

\definecolor{bottlebrush}{RGB}{211,39,55} % #d32737
\definecolor{alachua}{RGB}{242,169,0} % #F2A900
\definecolor{gator}{RGB}{34,136,76} % #22884C
\definecolor{darkblue}{RGB}{0,38,87} % #002657
\definecolor{perennial}{RGB}{106,42,96} % #6A2A60

\definecolor{bottlebrushLite}{RGB}{250,233,235} % #FAE9EB
\definecolor{alachuaLite}{RGB}{253,246,229} % #FDF6E5
\definecolor{gatorLite}{RGB}{232,243,237} % #E8F3ED
\definecolor{darkblueLite}{RGB}{229,233,238} % #E5E9EE
\definecolor{perennialLite}{RGB}{240,233,239} % #F0E9EF

\title{Pattern bounds for principal specializations of $\beta$-Grothendieck Polynomials}
\author{Hugh Dennin}
\date{June 2022}

\begin{document}

\maketitle

\begin{abstract}
    There has been recent interest in lower bounds for the principal specializations of Schubert polynomials $\nu_w := \fk S_w(1,\dots,1)$.
    We prove a conjecture of Yibo Gao in the setting of $1243$-avoiding permutations that gives a lower bound for $\nu_w$ in terms of the permutation patterns contained in $w$.
    We extended this result to principal specializations of $\beta$-Grothendieck polynomials $\nu^{(\beta)}_w := \fk G^{(\beta)}_w(1,\dots,1)$ by restricting to the class of vexillary $1243$-avoiding permutations.
    Our methods are bijective, offering a combinatorial interpretation of the coefficients $c_w$ and $c^{(\beta)}_w$ appearing in these conjectures.
\end{abstract}

\section{Introduction}
\emph{Schubert polynomials} were introduced by Lascoux and Sch\"utzenberger \cite{LS82a} as canonical representatives for Schubert cycles in the cohomology ring of the complete flag variety.
In particular, given a permutation $w\in S_n$ the Schubert polynomial $\fk S_w(x_1,\dots,x_{n-1})$ represents the cohomology class of the Schubert variety $X_w$.
% The \emph{principal specialization} $\fk S_w(1,\dots,1)$ is known to give the degree of the variety $X$

There has been recent interest in establishing bounds for the \emph{principal specialization} $\nu_w := \fk S_w(1,\dots,1) \in \Z_{\geq0}$, which is known to give the degree of the matrix Schubert variety corresponding to $w$ \cite{KM05}.
The first result in this vein came from Weigandt \cite{Wei18} who showed that $\nu_w \geq 1 + p_{132}(w)$, where, in general, $p_u(w)$ represents the number of occurrences of the pattern $u$ in $w$ (i.e., the number of subwords of $w$ whose entries appear in the same relative order as $u$).
This gave a conjecture of Stanley \cite{Sta17}, who theorized that $\nu_w = 2$ if and only if $p_{132}(w)=1$.

More recently \cite{Gao21}, Yibo Gao improved this lower bound by also including the number of occurrences of $1432$ patterns in $w$, that is, $\nu_w \geq 1 + p_{132}(w) + p_{1432}(w)$.
In the same paper, Gao conjectured that there should be a way to (in some sense) include every permutation pattern in a lower bound for $\nu_w$.
Define the following sequence of coefficients $c_w$ indexed by permutations $w$ recursively:
\begin{enumerate}
    \item Set $c_\emptyset := 1$ (here $\emptyset$ is the unique permutation of size $0$).
    \item Assuming that $c_u$ has been defined for all permutations $u$ of size less than $n$, for $w\in S_n$ set \[
        c_w := \nu_w - \sum_{u< w} c_up_u(w)
    \] where this sum is taken over all permutation patterns $u$ properly contained in $w$ (so $u\neq w$).
\end{enumerate}

\begin{conjecture}[Gao\footnote{
    Gao remarks that Conjecture \ref{conj:gao} was also observed independently by Christian Gaetz.
} \cite{Gao21}]
\label{conj:gao}
    For all permutations $w$, $c_w\geq 0$.
\end{conjecture}

Note that if Conjecture $\ref{conj:gao}$ were true, then we would get a lower bound \[
    \nu_w \geq \sum_{u<w} c_up_u(w).
\]
It can be checked that $c_\emptyset = c_{132} = c_{1432} = 1$ are the only nonzero values of $c_w$ for permutations $w$ of size $n\leq 4$, so this would extend the previous bounds of Weigandt and Gao.

Conjecture \ref{conj:gao} has been checked empirically to hold for all permutations $w\in S_n$ for $n\leq 8$ \cite{Gao21}.
As more evidence, M\'esz\'aros and Tanjaya \cite{MT21} showed that Conjecture \ref{conj:gao} is true whenever $w$ does not contain either of the patterns $1432$ or $1423$.

In this paper, we prove Conjecture \ref{conj:gao} for permutations $w$ that do not contain the pattern $1243$.
\begin{theorem}
\label{thm:gao_1243}
    For all $1243$-avoiding permutations $w$,
    $c_w\geq 0$.
\end{theorem}
% In Section \ref{sect:polys},
To illustrate a brief outline of the proof, it is known that $\nu_w$ enumerates a certain collection of diagrams corresponding to $w$ called \emph{reduced bumpless pipe dreams}.
To each of $w$'s reduced BPDs, we assign both a subword $v$ of $w$ and another reduced BPD corresponding to the permutation pattern $u$ that $v$ is an occurrence of (in $w$).
Fixing $v$,
this gives an injection from the reduced BPDs of $w$ assigned to $v$ into a certain subset of \emph{minimal} reduced BPDs of $u$.
This is bijective when $w$ avoids the pattern $1243$, in which case $c_w$ will enumerate the set of minimal reduced bumpless pipe dreams for $w$ and Theorem \ref{thm:gao_1243} will follow.
In general, this yields an upper bound for the principal specialization $\nu_w$ in terms the minimal reduced BPDs of its permutations patterns (Corollary \ref{cor:schub_ub}).
% We will see that this upper bound is exact in the $1243$-avoiding case, from which Theorem \ref{thm:gao_1243} will follow.

In parallel, we will consider the \emph{$\beta$-Grothendieck polynomials} $\fk G^{(\beta)}_w(x_1,\dots,x_{n-1})$ of Fomin and Kirillov \cite{FK94},
which specialize to Schubert polynomials in the case $\beta = 0$.
$\beta$-Grothendieck polynomials are known \cite{Hud14} to represent Schubert classes in the connective $K$-theory of the complete flag variety.
We can set up the same machinery for principal specializations of $\beta$-Grothendieck polynomials $\nu^{(\beta)}_w := \fk G^{(\beta)}_w(1,\dots,1)$,
which in this case are polynomials in $\Z_{\geq 0}[\beta]$,
by defining coefficients $c^{(\beta)}_w\in \Z[\beta]$ recursively as before:
\begin{enumerate}
    \item Set $c^{(\beta)}_\emptyset := 1$.
    \item Assuming that $c^{(\beta)}_u$ has been defined for all permutations $u$ of size less than $n$, for $w\in S_n$ set \[
        c^{(\beta)}_w :=
        \nu^{(\beta)}_w - \sum_{u<w} c^{(\beta)}_u p_u(w).
    \]
\end{enumerate}

Recall that a permutation is \emph{vexillary} if it avoids the pattern $2143$.
The techniques used to prove Theorem \ref{thm:gao_1243} are also applicable in the Grothendieck setting by restricting our view to a smaller pattern-avoidance class of permutations.
% to the case $\beta = 1$ (rather than $\beta = 0$) 
\begin{theorem}
\label{thm:groth_gao_1243_2143}
    For all vexillary $1243$-avoiding permutations $w$,
    $c^{(\beta)}_w \in \Z_{\geq0}[\beta]$.
\end{theorem}

The coefficients of $c^{(\beta)}_w$ have been confirmed by computer to be nonnegative for all permutations $w\in S_n$ of size $n\leq 9$.
This, coupled with the establishment of Theorem \ref{thm:groth_gao_1243_2143},
inspires us to make the following analogue of Conjecture \ref{conj:gao}.
\begin{conjecture}
\label{conj:groth_gao}
    For all permutations $w$, $c^{(\beta)}_w \in \Z_{\geq0}[\beta]$.
\end{conjecture}
We remark that $\nu_w$ and $c_w$ are just the constant terms of $\nu^{(\beta)}_w$ and $c^{(\beta)}_w$, respectively, so Conjecture \ref{conj:groth_gao} is a strengthening of Conjecture \ref{conj:gao}.
\section{Background \& Definitions}

% - - - - - - - - - - - - - - %
% - - PATTERN CONTAINMENT - - %
% - - - - - - - - - - - - - - %
\subsection{Pattern Containment}
\label{sect:pat}

Let $S_n$ denote the set of permutations of size $n$ i.e., the set of bijections from $[n] := \{1,2,\dots,n\}$ to itself.
We will write permutations as words in one-line notation.
In other words, for $w\in S_n$ a permutation, we identify $w$ with the word $w(1)w(2)\dots w(n)$.
Write $\abs{w} := n$ for the size of $w$.

By a \emph{subword} $v$ of a permutation $w\in S_n$, we mean a map $v: [m]\ra [n]$ such that there exist indices $1\leq s_1<\cdots<s_m\leq n$ with $v(i) = w(s_i)$ for all $1\leq i\leq m$.
As with permutations, we will write $v = v(1)v(2)\cdots v(m)$ and $\abs{v} := m$ for the size of $v$.
The relation $v\sseq w$ denotes the fact that $v$ is a subword of $w$ whereas $v\subset w$ denotes the same with the additional constraint $v\neq w$.
This notation suggests that we should think of $v$ as a subset of the indices of $w$.

%We first discuss the notations of permutation pattern containment and avoidance.
Let $w\in S_n$ be a permutation and let 
$v\sseq w$ be a subword of size $m$.
Define $\perm(v)$ to be the permutation $u\in S_m$ whose entries appear in the same relative order as the entries of $v$ i.e., $v(i) < v(j)$ if and only if $u(i) < u(j)$ for all $1\leq i<j\leq m$.
As an example, the subword $253$ of $12453$ has $\perm(253) = 132$ since $2<3<5$.

Conversely, given permutations $w\in S_n$ and $u\in S_m$,
we define $p_u(w)$ to be the number of subwords $v\sseq w$ satisfying $\perm(v) = u$.
% This quantity is known as the \emph{packing density} of $u$ in $w$.

\begin{definition}
\label{def:patterns}
    Let $w\in S_n$ and $u\in S_m$ be permutations.
    Say that $w$ \emph{avoids the pattern} $u$ or is \emph{$u$-avoiding} if $p_u(w) = 0$.
    Otherwise say that $w$ \emph{contains the pattern} $u$.
\end{definition}

Going back to our previous example, $21453$ contains the pattern $132$ (note that $132$ is not a \emph{subword} of $21453$).
On the other hand, $21453$ avoids the pattern $2143$.
Permutations that are $2143$-avoiding have a special name; they are known as \emph{vexillary} permutations.

We will use the notation $u\leq w$ to denote that $w$ contains the pattern $u$ and the notation $u< w$ to additionally indicate that $u\neq w$.\footnote{
    This is not to be confused with the \emph{Bruhat order} on $S_n$, which is also often denoted by $\leq$.
}
It is important to make the distinction between the partial orders $\sseq$, which is defined on \emph{subwords} of $w$, and $\leq$, which is defined on \emph{permutations} that appear as a pattern in $w$.

By convention, we define $[0]$ to be empty and let $\emptyset$ denote the unique permutation in $S_0$.
Every permutation $w$ contains $\emptyset$ as a pattern precisely once, that is $p_\emptyset(w) = 1$, and the occurrence of this pattern will also be denoted $\emptyset$.

% - - - - - - - - - - - - - - - %
% - - BUMPLESS PIPE DREAMS  - - %
% - - - - - - - - - - - - - - - %
\subsection{Bumpless Pipe Dreams}
\label{sect:bpd}

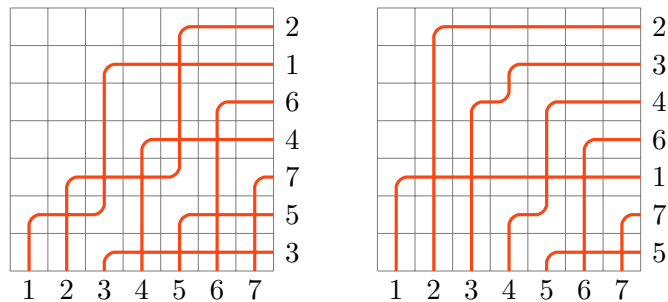
\begin{figure}[h]
    \centering
    \begin{tikzpicture}[scale=0.5]
    \draw[step=1.0,black,opacity=0.5,thin,xshift=0.5cm,yshift=0.5cm] (0,0) grid (7,7);
    \draw[\pipecolor,very thick]
    (1,0.5) to (1,1.75) to[bend left] (1.25,2) to (2.75,2) to[bend right] (3,2.25) to (3,5.75) to[bend left] (3.25,6) to (7.5,6) % pipe 1
    (2,0.5) to (2,2.75) to[bend left] (2.25,3) to (4.75,3) to[bend right] (5,3.25) to (5,6.75) to[bend left] (5.25,7) to (7.5,7) % pipe 2
    (3,0.5) to (3,0.75) to[bend left] (3.25,1) to (7.5,1) % pipe 3
    (4,0.5) to (4,3.75) to[bend left] (4.25,4) to (7.5,4) % pipe 4
    (5,0.5) to (5,1.75) to[bend left] (5.25,2) to (7.5,2) % pipe 5
    (6,0.5) to (6,4.75) to[bend left] (6.25,5) to (7.5,5) % pipe 6
    (7,0.5) to (7,2.75) to[bend left] (7.25,3) to (7.5,3); % pipe 7
    \foreach \i in {1,...,7}{\node at (\i,0) {\i};} % south labels
    \foreach \i\j in {7/2,6/1,5/6,4/4,3/7,2/5,1/3}{\node at (8,\i) {\j};} % east labels
    \end{tikzpicture}
    \qquad
    \begin{tikzpicture}[scale=0.5]
    \draw[step=1.0,black,opacity=0.5,thin,xshift=0.5cm,yshift=0.5cm] (0,0) grid (7,7);
    \draw[\pipecolor,very thick]
    (1,0.5) to (1,2.75) to[bend left] (1.25,3) to (7.5,3) % pipe 1
    (2,0.5) to (2,6.75) to[bend left] (2.25,7) to (7.5,7) % pipe 2
    (3,0.5) to (3,4.75) to[bend left] (3.25,5) to (3.75,5) to[bend right] (4,5.25) to (4,5.75) to[bend left] (4.25,6) to (7.5,6) % pipe 3
    (4,0.5) to (4,1.75) to[bend left] (4.25,2) to (4.75,2) to[bend right] (5,2.25) to (5,4.75) to[bend left] (5.25,5) to (7.5,5) % pipe 4
    (5,0.5) to (5,0.75) to[bend left] (5.25,1) to (7.5,1) % pipe 5
    (6,0.5) to (6,3.75) to[bend left] (6.25,4) to (7.5,4) % pipe 6
    (7,0.5) to (7,1.75) to[bend left] (7.25,2) to (7.5,2); % pipe 7
    \foreach \i in {1,...,7}{\node at (\i,0) {\i};} % south labels
    \foreach \i\j in {7/2,6/3,5/4,4/6,3/1,2/7,1/5}{\node at (8,\i) {\j};} % east labels
    \end{tikzpicture}
    \caption{Two bumpless pipe dreams of size $7$.}
    \label{fig:bpd}
\end{figure}

A $\emph{bumpless pipe dream}$ of size $n$ is a filling of the $n\times n$ grid using the tiles \[
    \bigblank\quad
    \bighorizontal\quad
    \bigvertical\quad
    \bigcross\quad
    \bigrelbow\quad
    \bigjelbow
\] such that the resulting network of ``pipes'' consists of $n$ pipes, each entering from a column along the south edge and exiting through a row along the east edge \cite{LLS21}.
In order from left-to-right, we will refer to these tiles as the \emph{blank}, \emph{dash}, \emph{bar}, \emph{cross}, \emph{r-elbow}, and \emph{j-elbow}, respectively.
See Figure \ref{fig:bpd} for two examples of bumpless pipe dreams of size $7$.

We let $\BPD(n)$ denote the collection of bumpless pipe dreams of size $n$.
The standard convention is to plot each $B\in \BPD(n)$ using matrix coordinates.
By $B_{i,j}$ we mean the tile in $B$ that is at position $(i,j)$.
The notation $B^{-1}(\blank)$ means the set of positions $(i,j)$ such that $B_{i,j} = \blank$
(and we make analogous definitions for $\horizontal$, $\vertical$, $\cross$, $\relbow$, and $\jelbow$).
Note that the pipes in a bumpless pipe dreams travel weakly in the northeast direction: a fact that we will make good use of in the sections to come.

We now make a brief digression to discuss a connection between bumpless pipe dreams and alternating sign matrices.
An \emph{alternating sign matrix} of size $n$ is a square $n\times n$ matrix $A$ consisting of the entries $0$, $+1$, and $-1$ such that in each row and column of $A$, the entries sum to $+1$ with the nonzero entries alternating in sign.
These conditions imply that each row and column of $A$ has at least one nonzero entry with the initial and final nonzero entries always being $+1$'s.
Let $\ASM(n)$ denote the set of alternating sign matrices of size $n$.

There is a well-known \cite{Wei21} % Lemma 3.1
bijection $\Phi: \BPD(n) \ra \ASM(n)$ where for $B$ a bumpless pipe dream of size $n$, we get an $n\times n$ alternating sign matrix via the formula \[
    \Phi(B)_{i,j} := \begin{cases}
        +1 & B_{i,j} = $\relbow$, \\
        -1 & B_{i,j} = $\jelbow$, \\
        0 & \text{else.}
    \end{cases}
\]
This means that the elbows in each row and each column of $B$ are alternating between $\relbow$ and $\jelbow$ with the initial and final elbows always being $\relbow$'s.
%In particular, there is always one fewer $\jelbow$ than there are $\relbow$s in a given row or column.
It will often be helpful to have this reductionist view of bumpless pipe dreams.

\begin{figure}[h]
    \centering
    \[\begin{pmatrix}
     0& 0& 0& 0&+1& 0& 0\\
     0& 0&+1& 0& 0& 0& 0\\
     0& 0& 0& 0& 0&+1& 0\\
     0& 0& 0& 0& 0& 0& 0\\
     0&+1& 0& 0&-1& 0&+1\\
    +1& 0&-1& 0&+1& 0& 0\\
     0& 0&+1& 0& 0& 0& 0
    \end{pmatrix}\]
    \caption{The alternating sign matrix $\Phi(B)$ corresponding to the first bumpless pipe dream $B$ in Figure \ref{fig:bpd}.}
    \label{fig:asm}
\end{figure}

We now return to our introduction of bumpless pipe dreams.
For $B\in \BPD(n)$, we get a permutation $w_B\in S_n$ by setting $w_B(x) := y$ if and only if $y\ra x$ is a pipe in $B$ (that is, the pipe in $B$ entering from column $y$ exits in row $x$).
To describe this relationship, we say that the bumpless pipe dream $B$ has permutation $w_B$.
% \emph{Coxeter product}
The bumpless pipe dreams in Figure \ref{fig:bpd} have permutations $2164753$ and $2346175$, respectively.
By $\BPD(w)$ we mean the set of bumpless pipe dreams that have permutation $w$.

\begin{figure}[h]
    \centering
    \begin{tikzpicture}[scale=0.5]
    \draw[step=1.0,black,opacity=0.5,thin,xshift=0.5cm,yshift=0.5cm] (0,0) grid (7,7);
    \draw[UForange,very thick]
    (1,0.5) to (1,1.75) to[bend left] (1.25,2) to (2.75,2) to[bend right] (3,2.25) to (3,2.75) to[bend left] (3.25,3) to (4.75,3) to[bend right] (5,3.25) to (5,3.75) to[bend left] (5.25,4) to (7.5,4) % pipe 1
    (2,0.5) to (2,2.75) to[bend left] (2.25,3) to (2.75,3) to[bend right] (3,3.25) to (3,5.75) to[bend left] (3.25,6) to (7.5,6) % pipe 2
    (3,0.5) to (3,0.75) to[bend left] (3.25,1) to (7.5,1) % pipe 3
    (4,0.5) to (4,3.75) to[bend left] (4.25,4) to (4.75,4) to[bend right] (5,4.25) to (5,6.75) to[bend left] (5.25,7) to (7.5,7) % pipe 4
    (5,0.5) to (5,1.75) to[bend left] (5.25,2) to (7.5,2) % pipe 5
    (6,0.5) to (6,4.75) to[bend left] (6.25,5) to (7.5,5) % pipe 6
    (7,0.5) to (7,2.75) to[bend left] (7.25,3) to (7.5,3); % pipe 7
    \foreach \i in {1,...,7}{\node at (\i,0) {\i};} % south labels
    \foreach \i\j in {7/4,6/2,5/6,4/1,3/7,2/5,1/3}{\node at (8,\i) {\j};} % east labels
    \end{tikzpicture}
    \caption{The resulting network of pipes $B_K$ obtained from the first bumpless pipe dream $B$ in Figure \ref{fig:bpd}.}
    \label{fig:bpd_K}
\end{figure}
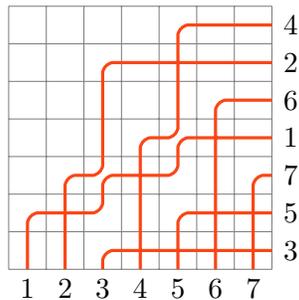

There is an alternative way to interpret how the network of pipes in $B$ behaves at crosses; namely, we consider two pipes to cross only at the first $\cross$ they partake in, otherwise diverging away from each other at subsequent crosses.
To make this precise, order $B^{-1}(\cross)$ from left-to-right, bottom-to-top.
Then for each $(i,j)$ in $B^{-1}(\cross)$ in order, check if the pipes passing through the cross $B_{i,j}$ have already crossed at a previous position in $B^{-1}(\cross)$.
If so, change $B_{i,j}$ from a cross to a \emph{bump}:\footnote{
    This eponymous tile, absent in the bumpless pipe dream model, is present in the older \emph{pipe dream} model.
} \[
    \bigcross\quad
    \raisebox{4pt}{$\longmapsto$}\quad
    \bigbump
\]
Let $B_K$ denote the diagram obtained by this process.\footnote{
    The $K$ refers to $K$-theory.
}
See Figure \ref{fig:bpd_K} for an example.

We get another permutation in $S_n$ associated to $B$, denoted $\partial(B)$, by setting $\partial(B)(x) := y$ if and only if $y\ra x$ is a pipe appearing in $B_K$.
Say that $B$ has \emph{type} $\partial(B)$.
% \emph{Demazure product}
Keeping with our example, the first bumpless pipe dream in Figure \ref{fig:bpd} has type $4261753$, as evidenced in Figure \ref{fig:bpd_K}.
The second BPD in Figure \ref{fig:bpd} has no double crossings between pipes and hence remains unchanged by this process, so it has type $2356175$.
Let $\BPD_K(w)$ denote the set of bumpless pipe dreams of type $w$.

% - - - - - - - - - - - - - - - - - %
% - - MINIMAL AND REDUCED BPDs  - - %
% - - - - - - - - - - - - - - - - - %
\subsection{Minimal and Reduced BPDs}
\label{sect:min_red}

\begin{definition}
    Let $B\in \BPD(n)$.
    Say that a pipe $y\ra x$ is \emph{removable} in $B$ if $B_{x,y} = \relbow$, and there are no other $\relbow$'s in row $x$ and column $y$.
    Say that $B$ is \emph{minimal} if it contains no removable pipes.
\end{definition}

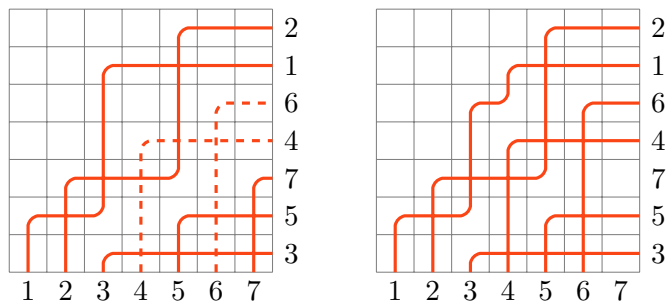
\begin{figure}[h]
    \centering
    \begin{tikzpicture}[scale=0.5]
    \draw[step=1.0,black,opacity=0.5,thin,xshift=0.5cm,yshift=0.5cm] (0,0) grid (7,7);
    \draw[\pipecolor,very thick]
    (1,0.5) to (1,1.75) to[bend left] (1.25,2) to (2.75,2) to[bend right] (3,2.25) to (3,5.75) to[bend left] (3.25,6) to (7.5,6) % pipe 1
    (2,0.5) to (2,2.75) to[bend left] (2.25,3) to (4.75,3) to[bend right] (5,3.25) to (5,6.75) to[bend left] (5.25,7) to (7.5,7) % pipe 2
    (3,0.5) to (3,0.75) to[bend left] (3.25,1) to (7.5,1) % pipe 3
    (5,0.5) to (5,1.75) to[bend left] (5.25,2) to (7.5,2) % pipe 5
    (7,0.5) to (7,2.75) to[bend left] (7.25,3) to (7.5,3); % pipe 7
    \draw[\pipecolor,very thick,dashed]
    (4,0.5) to (4,3.75) to[bend left] (4.25,4) to (7.5,4) % pipe 4
    (6,0.5) to (6,4.75) to[bend left] (6.25,5) to (7.5,5); % pipe 6
    \foreach \i in {1,...,7}{\node at (\i,0) {\i};} % south labels
    \foreach \i\j in {7/2,6/1,5/6,4/4,3/7,2/5,1/3}{\node at (8,\i) {\j};} % east labels
    \end{tikzpicture}
    \qquad
    \begin{tikzpicture}[scale=0.5]
    \draw[step=1.0,black,opacity=0.5,thin,xshift=0.5cm,yshift=0.5cm] (0,0) grid (7,7);
    \draw[\pipecolor,very thick]
    (1,0.5) to (1,1.75) to[bend left] (1.25,2) to (2.75,2) to[bend right] (3,2.25) to (3,4.75) to[bend left] (3.25,5) to (3.75,5) to[bend right] (4,5.25) to (4,5.75) to[bend left] (4.25,6) to (7.5,6) % pipe 1
    (2,0.5) to (2,2.75) to[bend left] (2.25,3) to (4.75,3) to[bend right] (5,3.25) to (5,6.75) to[bend left] (5.25,7) to (7.5,7) % pipe 2
    (3,0.5) to (3,0.75) to[bend left] (3.25,1) to (7.5,1) % pipe 3
    (4,0.5) to (4,3.75) to[bend left] (4.25,4) to (7.5,4) % pipe 4
    (5,0.5) to (5,1.75) to[bend left] (5.25,2) to (7.5,2) % pipe 5
    (6,0.5) to (6,4.75) to[bend left] (6.25,5) to (7.5,5); % pipe 6
    (7,0.5) to (7,2.75) to[bend left] (7.25,3) to (7.5,3); % pipe 7
    \foreach \i in {1,...,7}{\node at (\i,0) {\i};} % south labels
    \foreach \i\j in {7/2,6/1,5/6,4/4,3/7,2/5,1/3}{\node at (8,\i) {\j};} % east labels
    \end{tikzpicture}
    \caption{On the left, the first bumpless pipe dream from Figure \ref{fig:bpd} with its two removable pipes indicated via dashed lines. A slight modification has been made to this BPD on the right to make it minimal.}
    \label{fig:min_bpd}
\end{figure}

Observe that if pipe $y\ra x$ is removable in $B\in \BPD(n)$, then there are also not any $\jelbow$'s in row $x$ and column $y$ since the elbows in $B$ are alternating.
This shows that the $\relbow$ at $(x,y)$ belongs to pipe $y\ra x$, so pipe $y\ra x$ is \emph{undrooped} or hook-shaped (it does not contain a $\jelbow$).
An example is illustrated in Figure \ref{fig:min_bpd}.
We will have $\mBPD(n)$ denote the collection of minimal bumpless pipe dreams of size $n$.

In the language of alternating sign matrices, as a pipe $y\ra x$ is removable in $B$ if and only if $\Phi(B)_{x,y} = +1$, and this is the only nonzero entry in row $x$ and column $y$ of $\Phi(B)$.

A good heuristic to have is that removable pipes somehow do not contribute any complexity to a bumpless pipe dream.
This will be made explicit in Section \ref{sect:rem_pipes}, where we will investigate a process for removing the aptly named removable pipes from a BPD.

Fix $w\in S_n$.
For $B\in \BPD(w)$ with removable pipes $y_1\ra x_1,\dots,y_r\ra x_r$,
we get a subword $v\sseq w$ of size $n-r$ simply by ignoring the values $y_1,\dots,y_r$ in $w$.
For the BPD on the left in Figure \ref{fig:min_bpd}, we get the subword $21753$ of $2164753$.
Let $\BPD(w;v)$ denote the collection of all $B\in \BPD(w)$ that correspond to the subword $v\sseq w$ in this way.
This gives a partition \[
    \BPD(w) = \bigsqcup_{v\sseq w} \BPD(w;v).
\]
In particular, observe that $\BPD(w;w) = \mBPD(w)$ where this latter set is, of course, the collection of minimal bumpless pipe dreams that have permutation $w$.

%We let $\BPD_m(w)$ (resp. $\bpd_m(w)$) denote the set of minimal bumpless pipe dreams (resp. minimal reduced bumpless pipe dreams) corresponding to $w$.

%Formally, a bumpless pipe dream $B\in\BPD$ is a function $[n]\times[n]\ra \tiles$.
%It is well known that $\mathfrak S_w(1,\dots,1) = |\BPD(w)|$.

%We will think of each $B\in \BPD$ as a collection of pipes indexed by $[n]$, where by $p_j$ we mean the pipe entering $B$ from the bottom of column $j$ and exiting to the right in row $i:= w^{-1}(j)$.

The next definition is standard in the literature.

\begin{definition}
    We say that a bumpless pipe dream is \emph{reduced} if no two of its pipes cross more than once.
\end{definition}

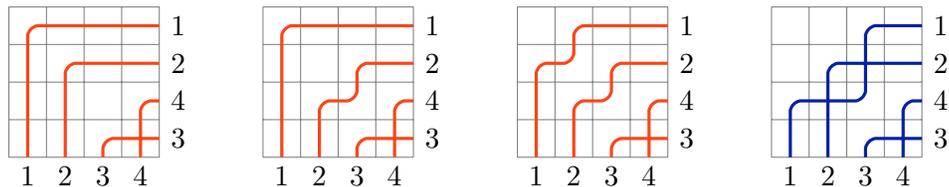
\begin{figure}[h]
    \centering
    \begin{tikzpicture}[scale=0.5]
    \draw[step=1.0,black,opacity=0.5,thin,xshift=0.5cm,yshift=0.5cm] (0,0) grid (4,4);
    \draw[\pipecolor,very thick]
    (1,0.5) to (1,3.75) to[bend left] (1.25,4) to (4.5,4) % pipe 1
    (2,0.5) to (2,2.75) to[bend left] (2.25,3) to (4.5,3) % pipe 2
    (3,0.5) to (3,0.75) to[bend left] (3.25,1) to (4.5,1) % pipe 3
    (4,0.5) to (4,1.75) to[bend left] (4.25,2) to (4.5,2); % pipe 4
    \foreach \i in {1,...,4}{\node at (\i,0) {\i};} % south labels
    \foreach \i\j in {4/1,3/2,2/4,1/3}{\node at (5,\i) {\j};} % east labels
    \end{tikzpicture}
    \qquad
    \begin{tikzpicture}[scale=0.5]
    \draw[step=1.0,black,opacity=0.5,thin,xshift=0.5cm,yshift=0.5cm] (0,0) grid (4,4);
    \draw[\pipecolor,very thick]
    (1,0.5) to (1,3.75) to[bend left] (1.25,4) to (4.5,4) % pipe 1
    (2,0.5) to (2,1.75) to[bend left] (2.25,2) to (2.75,2) to[bend right] (3,2.25) to (3,2.75) to[bend left] (3.25,3) to (4.5,3) % pipe 2
    (3,0.5) to (3,0.75) to[bend left] (3.25,1) to (4.5,1) % pipe 3
    (4,0.5) to (4,1.75) to[bend left] (4.25,2) to (4.5,2); % pipe 4
    \foreach \i in {1,...,4}{\node at (\i,0) {\i};} % south labels
    \foreach \i\j in {4/1,3/2,2/4,1/3}{\node at (5,\i) {\j};} % east labels
    \end{tikzpicture}
    \qquad
    \begin{tikzpicture}[scale=0.5]
    \draw[step=1.0,black,opacity=0.5,thin,xshift=0.5cm,yshift=0.5cm] (0,0) grid (4,4);
    \draw[\pipecolor,very thick]
    (1,0.5) to (1,2.75) to[bend left] (1.25,3) to (1.75,3) to[bend right] (2,3.25) to (2,3.75) to[bend left] (2.25,4) to (4.5,4) % pipe 1
    (2,0.5) to (2,1.75) to[bend left] (2.25,2) to (2.75,2) to[bend right] (3,2.25) to (3,2.75) to[bend left] (3.25,3) to (4.5,3) % pipe 2
    (3,0.5) to (3,0.75) to[bend left] (3.25,1) to (4.5,1) % pipe 3
    (4,0.5) to (4,1.75) to[bend left] (4.25,2) to (4.5,2); % pipe 4
    \foreach \i in {1,...,4}{\node at (\i,0) {\i};} % south labels
    \foreach \i\j in {4/1,3/2,2/4,1/3}{\node at (5,\i) {\j};} % east labels
    \end{tikzpicture}
    \qquad%\qquad
    \begin{tikzpicture}[scale=0.5]
    \draw[step=1.0,black,opacity=0.5,thin,xshift=0.5cm,yshift=0.5cm] (0,0) grid (4,4);
    \draw[UFblue,very thick]
    (1,0.5) to (1,1.75) to[bend left] (1.25,2) to (2.75,2) to[bend right] (3,2.25) to (3,3.75) to[bend left] (3.25,4) to (4.5,4) % pipe 1
    (2,0.5) to (2,2.75) to[bend left] (2.25,3) to (4.5,3) % pipe 2
    (3,0.5) to (3,0.75) to[bend left] (3.25,1) to (4.5,1) % pipe 3
    (4,0.5) to (4,1.75) to[bend left] (4.25,2) to (4.5,2); % pipe 4
    \foreach \i in {1,...,4}{\node at (\i,0) {\i};} % south labels
    \foreach \i\j in {4/1,3/2,2/4,1/3}{\node at (5,\i) {\j};} % east labels
    \end{tikzpicture}
    \caption{The four bumpless pipe dreams with permutation $1243$, the first three of which are reduced (in orange).}
    \label{fig:red_bpd}
\end{figure}

Note that the second BPD in Figure \ref{fig:bpd} is reduced while the first is not. In Figure \ref{fig:red_bpd}, we illustrate the reduced and nonreduced bumpless pipe dreams that have permutation $1243$.

Observe that for $B\in \BPD(n)$, $B_K$ differs from $B$ precisely when $B$ is not reduced i.e., there are two pipes in $B$ that cross more than once.
When $B$ is reduced, we get equality $B = B_K$ and so the permutation and type of $B$ coincide
(as we saw with the second BPD in Figure \ref{fig:bpd}).

For each definition using the notation $\BPD$, we make an analogous definition for reduced bumpless pipe dreams using the notation $\bpd$.
Explicitly, we define $\bpd(n)$ (resp. $\mbpd(n)$, $\bpd(w)$, $\bpd_K(w)$, $\mbpd(w)$, $\bpd(w;v)$) to be the subset of reduced bumpless pipe dreams in $\BPD(n)$ (resp. $\mBPD(n)$, $\BPD(w)$, $\BPD_K(w)$, $\mBPD(w)$, $\BPD(w;v)$).
Note that $\bpd(w) = \bpd_K(w)$ by the previous paragraph, so this definition is redundant.
As before, we have the equalities \[
    \bpd(w) = \bigsqcup_{v\sseq w} \bpd(w;v)
\] and $\bpd_w(w) = \mbpd(w)$.

To help organize these notions, we provide the following diagram of containments:
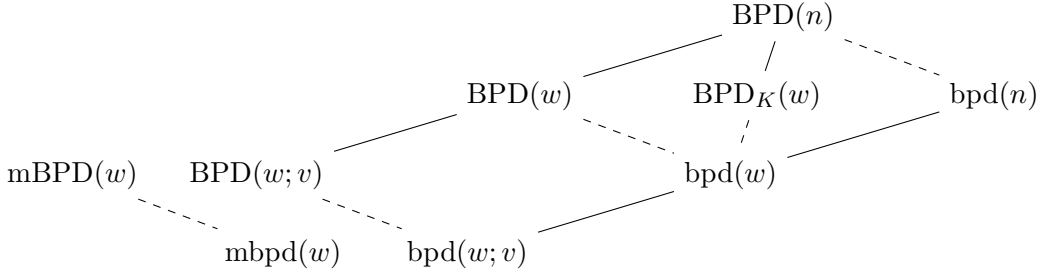
\[\begin{tikzpicture}[scale=.7]
    \node (A) at (10,4.5) {$\BPD(n)$};
    \node (B) at (5,3) {$\BPD(w)$};
    \node (B') at (9.5,3) {$\BPD_K(w)$};
    \node (C) at (0,1.5) {$\BPD(w;v)$};
    \node (C') at (-3.5,1.5) {$\mBPD(w)$};
    \node (a) at (14,3) {$\bpd(n)$};
    \node (b) at (9,1.5) {$\bpd(w)$};
    \node (c) at (4,0) {$\bpd(w;v)$};
    \node (c') at (0.5,0) {$\mbpd(w)$};
    \draw
    (A) -- (B) (A) -- (B')
    (B) -- (C) % (B) -- (C')
    (a) -- (b)
    (b) -- (c); % (b) -- (c');
    \draw[dashed]
    (A) -- (a)
    (B) -- (b)
    (B') -- (b)
    (C) -- (c)
    (C') -- (c');
\end{tikzpicture}\]

% - - - - - - - - - - - - - - - - - - - - - - - %
% - - GROTHENDIECK AND SCHUBERT POLYNOMIALS - - %
% - - - - - - - - - - - - - - - - - - - - - - - %
\subsection{Grothendieck and Schubert Polynomials}
\label{sect:polys}

Throughout this subsection, fix a permutation $w\in S_n$.

The \emph{(single) $\beta$-Grothendieck polynomial} for $w$, as introduced by Fomin and Kirillov \cite{FK94},
is a polynomial in $\Z_{\geq0}[x_1,x_2,\dots][\beta]$ known \cite{Las02,Wei21} % \cite{Las02}
to be given by the formula \begin{equation}
\label{eq:groth}
    \fk G^{(\beta)}_w(x_1,\dots,x_{n-1}) =
    \beta^{-\ell(w)} \sum_{B\in \BPD_K(w)}
    \Bigg(\prod_{(i,j)\in B^{-1}(\smallblank)} \beta x_i
    \prod_{(i,j)\in B^{-1}(\smalljelbow)} (1+\beta x_i) \Bigg).
\end{equation}
Here $\ell(w)$ denotes the \emph{Coxeter length} of $w$ which is the number of \emph{inversions} in $w$ i.e., the number of tuples $(i,j)$ where $1\leq i<j\leq n$ and $w(i)>w(j)$.
We remark that $\ell(w)$ is also equal to the pattern count $p_{21}(w)$.
It is well known that $\abs{B^{-1}(\blank)} \geq \ell(w)$ with equality if and only if $B$ is reduced \cite{Wei21}, % Lemma 3.7
so $\fk G^{(\beta)}_w(x_1,\dots,x_{n-1})$ is indeed a polynomial in $\Z_{\geq 0}[x_1,x_2,\dots][\beta]$.

Of interest to us is the \emph{principal specialization} of the $\beta$-Grothendieck polynomial $\nu_w^{(\beta)} := \fk G_w^{(\beta)}(1,\dots,1)$, which itself is a polynomial in $\Z_{\geq 0}[\beta]$.
From the above, $\nu^{(\beta)}_w$ can be calculated as \[
    \nu^{(\beta)}_w = \fk G^{(\beta)}_w(1,\dots,1) =
    \beta^{-\ell(w)} \sum_{B\in \BPD_K(w)}
    \beta^{|B^{-1}(\smallblank)|}
    (1+\beta)^{|B^{-1}(\smalljelbow)|}.
\]
To simplify notation, define the \emph{$\beta$-weight} of a bumpless pipe dream $B\in \BPD(n)$ to be the quantity $\wt^{(\beta)}(B) := \beta^{|B^{-1}(\smallblank)|-\ell(w)} (1+\beta)^{|B^{-1}(\smalljelbow)|} \in \Z_{\geq 0}[\beta]$.
Extend this weighting to collections of bumpless pipe dreams $\mathcal S\sseq \BPD(n)$ by setting $\wt^{(\beta)}(\mathcal S) := \sum_{B\in \mathcal S} \wt^{(\beta)}(B)$.
Then $\nu^{(\beta)}_w = \wt^{(\beta)}(\BPD_K(w))$.

% that is, the evaluation $\nu_w^{(\beta)} := \fk G_w^{(\beta)}(1,\dots,1)$.
% where $x_i=1$ for all $1\leq i\leq n-1$.
% Let $\nu_w^{(\beta)} := \fk G_w^{(\beta)}(1,\dots,1)$ denote the principle specialization of the Grothendieck polynomial of $w$ for the parameter $\beta$.

There are two values of $\beta$ that will be of special interest to us, the first being the evaluation $\beta := 0$.
Let $B\in \BPD_K(w)$ be a bumpless pipe dream of type $w$.
The monomial contributed by $B$ in Equation \ref{eq:groth} has a factor of $\beta$ if and only if $\abs{B^{-1}(\blank)} > \ell(w)$, which, as we saw, occurs precisely when $B$ is not reduced.
Thus setting $\beta := 0$ filters our sum for $\fk G^{(\beta)}_w(x_1,\dots,x_{n-1})$ as \[
    \fk S_w(x_1,\dots,x_{n-1}) :=
    \fk G^{(0)}_w(x_1,\dots,x_{n-1}) =
    \sum_{B\in \bpd(w)}
    \Bigg(\prod_{(i,j)\in B^{-1}(\smallblank)} x_i \Bigg).
\]
This is a known \cite{LLS21} formula for the \emph{(single) Schubert polynomial} corresponding to $w$ of Lascoux and Sch\"utzenberger \cite{LS82a}.
Writing $\nu_w := \nu^{(0)}_w$, we see that the principal specialization of $\fk S_w(x_1,\dots,x_{n-1})$ gives an enumeration \[
    \nu_w = \fk S_w(1,\dots,1) = \abs{\bpd(w)}.
\]

% In a similar vein, the evaluation $\beta=-1$ recovers the regular \emph{(single) Grothendieck polynomial} for $w$, as introduced by Lascoux and Sch\"utzenberger \cite{LS82b}.

The second specialization of $\beta$ that we will examine is the choice $\beta := 1$.
% For $B\in \BPD(w)$, let $J(B)$ denote the quantity $\abs{B^{-1}(\jelbow)}$ i.e., the number of j-elbows in $B$.
Observe that the $1$-weight of a bumpless pipe dream $B\in \BPD(n)$ is simply $\wt^{(1)}(B) = 2^{|B^{-1}(\smalljelbow)|}$.
This weight is interesting when considering the bijection to alternating sign matrices where,
since the number of $\jelbow$'s in $B$ is the same as the number of $-1$'s in $\Phi(B)$, $\wt^{(1)}(B)$ is the same as the weight of $\Phi(B)$ in what is known as the \emph{$2$-enumeration} of $\ASM(n)$. % \cite{?}

It should be noted that the $\beta$-Grothendieck principal specializations $\nu^{(\beta)}_w$ have previously been studied for $w$ vexillary, where a determinantal upper bound was established in terms of Delannoy paths on Young tableaux \cite{MPP22}.
There is a $q$-analog of $\nu^{(\beta)}_w$ given by $\fk G^{(\beta)}_w(1,q,\dots,q^{n-2})$ where we instead send $x_i\mapsto q^{i-1}$.
Fomin and Stanley \cite{FS94} proved an elegant formula (originally conjectured by Macdonald \cite{Mac91}) for the Schubert case $\fk S_w(1,q,\dots,q^{n-2})$ as a weighted sum over the reduced words for $w$.
This specializes to a formula for $\nu_w$ by taking $q=1$.
Marberg and Pawlowski \cite{MP21} proved a analogous formula for the $x_i\mapsto q^{i-1}$ specializations of \emph{back stable} Schubert and $\beta$-Grothendieck polynomials, which are further generalizations of ordinary Schubert and $\beta$-Grothendieck polynomials that admit variables $x_i$ with indices $i\leq 0$.

% \[
%     \nu^{(1)}_w = \fk S^{(1)}_w(1,\dots,1)
%     = \sum_{B\in\BPD_K(w)} 2^{J(B)}
% \] where the weight $2^{J(B)}$ for $B\in \BPD_K(w)$ contributes to what is known as the \emph{$2$-enumeration} of $\BPD_K(w)$.
% Use the notation $\norm{\mathcal S}_2$ to denote the $2$-enumeration of a subset $\mathcal S\sseq \BPD(n)$ i.e., \[
%     \norm{\mathcal S}_2 := \sum_{B\in \mathcal S} 2^{J(B)}.
% \]

% and consider pipes $y_1\ra x_1$ and $y_2\ra x_2$ in $B$ where $x_1< x_2$.
% These pipes cross an odd number of times if and only if $y_1> y_2$ i.e., if and only if $(x_1,x_2)$ is an inversion in $w$.
% It follows that $\abs{B^{-1}(\cross)}\geq \ell(w)$.
\section{Results}

% - - - - - - - - - - - - %
% - - REMOVING PIPES  - - %
% - - - - - - - - - - - - %
\subsection{Removing Pipes}
\label{sect:rem_pipes}

Fix a permutation $w\in S_n$ and a subword $v\leq w$ of size $m$, and set $u := \perm(v) \in S_m$.
Write $y_1\ra x_1,\dots,y_r\ra x_r$ for the removable pipes in any bumpless pipe dream in $\BPD(w;v)$
(so $r := n-m$ and $y_1,\dots,y_r$ are the entries in $w$ not appearing in $v$).
Assume $y_1< \cdots< y_r$.

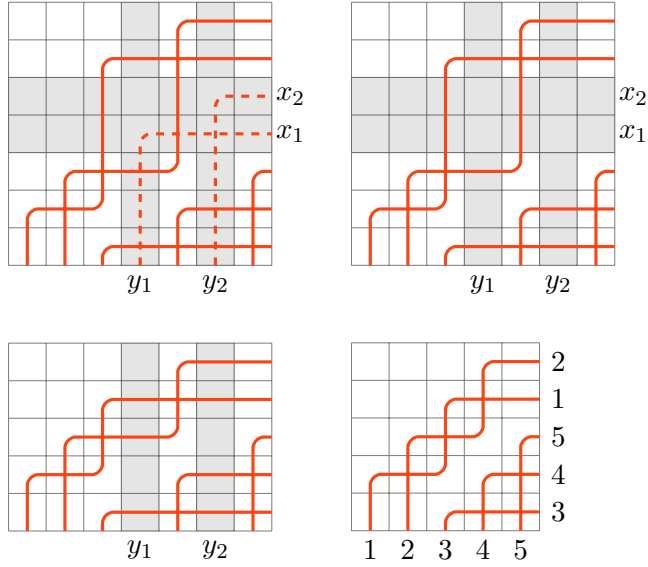
\begin{figure}[h]
    \centering
    \begin{tabular}{l l}
    \begin{tikzpicture}[scale=0.5]
    \filldraw[gray,opacity=0.2]
    (0.5,3.5) rectangle (7.5,5.5)
    (3.5,0.5) rectangle (4.5,7.5)
    (5.5,0.5) rectangle (6.5,7.5);
    \draw[step=1.0,black,opacity=0.5,thin,xshift=0.5cm,yshift=0.5cm] (0,0) grid (7,7);
    \draw[\pipecolor,very thick]
    (1,0.5) to (1,1.75) to[bend left] (1.25,2) to (2.75,2) to[bend right] (3,2.25) to (3,5.75) to[bend left] (3.25,6) to (7.5,6) % pipe 1
    (2,0.5) to (2,2.75) to[bend left] (2.25,3) to (4.75,3) to[bend right] (5,3.25) to (5,6.75) to[bend left] (5.25,7) to (7.5,7) % pipe 2
    (3,0.5) to (3,0.75) to[bend left] (3.25,1) to (7.5,1) % pipe 3
    (5,0.5) to (5,1.75) to[bend left] (5.25,2) to (7.5,2) % pipe 5
    (7,0.5) to (7,2.75) to[bend left] (7.25,3) to (7.5,3); % pipe 7
    \draw[\pipecolor,very thick,dashed]
    (4,0.5) to (4,3.75) to[bend left] (4.25,4) to (7.5,4) % pipe 4
    (6,0.5) to (6,4.75) to[bend left] (6.25,5) to (7.5,5); % pipe 6
    \node at (4,0) {$y_1$};
    \node at (6,0) {$y_2$};
    \node at (8,4) {$x_1$};
    \node at (8,5) {$x_2$};
    %\foreach \i in {1,...,7}{\node at (\i,0) {\i};} % south labels
    %\foreach \i\j in {7/2,6/1,5/6,4/4,3/7,2/5,1/3}{\node at (8,\i) {\j};} % east labels
    \end{tikzpicture}
    \vspace{1em}
    &
    \begin{tikzpicture}[scale=0.5]
    \filldraw[gray,opacity=0.2]
    (0.5,3.5) rectangle (7.5,5.5)
    (3.5,0.5) rectangle (4.5,7.5)
    (5.5,0.5) rectangle (6.5,7.5);
    \draw[black,opacity=0.5,thin,xshift=0.5cm,yshift=0.5cm] (0,0) grid (7,7);
    \draw[\pipecolor,very thick]
    (1,0.5) to (1,1.75) to[bend left] (1.25,2) to (2.75,2) to[bend right] (3,2.25) to (3,5.75) to[bend left] (3.25,6) to (7.5,6) % pipe 1
    (2,0.5) to (2,2.75) to[bend left] (2.25,3) to (4.75,3) to[bend right] (5,3.25) to (5,6.75) to[bend left] (5.25,7) to (7.5,7) % pipe 2
    (3,0.5) to (3,0.75) to[bend left] (3.25,1) to (7.5,1) % pipe 3
    (5,0.5) to (5,1.75) to[bend left] (5.25,2) to (7.5,2) % pipe 5
    (7,0.5) to (7,2.75) to[bend left] (7.25,3) to (7.5,3); % pipe 7
    \node at (4,0) {$y_1$};
    \node at (6,0) {$y_2$};
    \node at (8,4) {$x_1$};
    \node at (8,5) {$x_2$};
    \end{tikzpicture}
    \\
    \begin{tikzpicture}[scale=0.5]
    \filldraw[gray,opacity=0.2]
    (3.5,0.5) rectangle (4.5,5.5)
    (5.5,0.5) rectangle (6.5,5.5);
    \draw[black,opacity=0.5,thin,xshift=0.5cm,yshift=0.5cm] (0,0) grid (7,5);
    \draw[\pipecolor,very thick]
    (1,0.5) to (1,1.75) to[bend left] (1.25,2) to (2.75,2) to[bend right] (3,2.25) to (3,3.75) to[bend left] (3.25,4) to (7.5,4) % pipe 1
    (2,0.5) to (2,2.75) to[bend left] (2.25,3) to (4.75,3) to[bend right] (5,3.25) to (5,4.75) to[bend left] (5.25,5) to (7.5,5) % pipe 2
    (3,0.5) to (3,0.75) to[bend left] (3.25,1) to (7.5,1) % pipe 3
    (5,0.5) to (5,1.75) to[bend left] (5.25,2) to (7.5,2) % pipe 5
    (7,0.5) to (7,2.75) to[bend left] (7.25,3) to (7.5,3); % pipe 7
    \node at (4,0) {$y_1$};
    \node at (6,0) {$y_2$};
    \end{tikzpicture}
    &
    \begin{tikzpicture}[scale=0.5]
    \draw[step=1.0,black,opacity=0.5,thin,xshift=0.5cm,yshift=0.5cm] (0,0) grid (5,5);
    \draw[\pipecolor,very thick]
    (1,0.5) to (1,1.75) to[bend left] (1.25,2) to (2.75,2) to[bend right] (3,2.25) to (3,3.75) to[bend left] (3.25,4) to (5.5,4) % pipe 1
    (2,0.5) to (2,2.75) to[bend left] (2.25,3) to (3.75,3) to[bend right] (4,3.25) to (4,4.75) to[bend left] (4.25,5) to (5.5,5) % pipe 2
    (3,0.5) to (3,0.75) to[bend left] (3.25,1) to (5.5,1) % pipe 3
    (4,0.5) to (4,1.75) to[bend left] (4.25,2) to (5.5,2) % pipe 5
    (5,0.5) to (5,2.75) to[bend left] (5.25,3) to (5.5,3); % pipe 7
    \foreach \i\j in {1,...,5}{\node at (\i,0) {\i};} % south labels
    \foreach \i\j in {5/2,4/1,3/5,2/4,1/3}{\node at (6,\i) {\j};} % east labels
    \end{tikzpicture}
    \end{tabular}
    \caption{The map $\varphi_v^w$ applied to the first bumpless pipe dream $B$ in Figure \ref{fig:bpd} ($w = 2164753$, $v = 21753$). The final diagram $\varphi_v^w(B)$ is a minimal bumpless pipe dream with permutation $u := \perm(v) = 21543$.}
    \label{fig:phi_v^w}
\end{figure}

Let $B\in \BPD(w;v)$.
We now describe a procedure for \emph{removing} the (suggestively named) removable pipes from $B$:
\begin{enumerate}
    \item Erase the pipes $y_1\ra x_1,\dots,y_r\ra x_r$ from $B$.
    % (without disturbing any non-removable pipes that cross any of these pipes).
    \item Since each pipe $y_k\ra x_k$ ($1\leq k\leq r$) was removable, the remaining tiles in row $x_k$ consist entirely of $\blank$ and $\vertical$ tiles.
    We can now \emph{contract} the rows $x_1,\dots,x_r$ whilst maintaining the connectivity of the pipes that pass through them
    (i.e., delete the cells in these rows from $B$ and then top-justify the resulting diagram).
    \item Likewise, the remaining tiles in each column $y_k$ ($1\leq k\leq r$) consist entirely of $\blank$ and $\horizontal$ tiles.
    We can now \emph{contract} the columns $y_1,\dots,y_r$ whilst maintaining the connectivity of the pipes that pass through them
    (i.e., delete the cells in these columns from $B$ and then left-justify the resulting diagram).
\end{enumerate}
Let $\varphi_v^w(B)$ denote the diagram produced by the above procedure.
See Figure \ref{fig:phi_v^w} for an example.

It should be noted that $\varphi_v^w(B)$ could just as easily have been obtained by taking the $m\times m$ ``submatrix'' of $B$ by deleting the rows $x_1,\dots,x_r$ and columns $y_1,\dots,y_r$ from $B$ (as if taking a $m\times m$ minor of a matrix).
Indeed, in the language of alternating sign matrices $\Phi(\varphi_v^w(B))$ is literally this $m\times m$ submatrix of $\Phi(B)$.

One benefit of stratifying this procedure as we have done is that it becomes clear why $\varphi_v^w(B)$ is a bumpless pipe dream of size $m$ that, in particular, has permutation $u$.
This first statement immediate since all the pipes remain connected.
To see the latter, let $s_1<\cdots<s_m$ denote the indices of $v$ in $w$ and let $t_1<\cdots<t_m$ denote the entries of $v$.
Then $\varphi_v^w(B)_{i,j} = B_{s_i,t_j}$ for all $1\leq i,j\leq m$.
Furthermore, if $\varphi_v^w(B)$ has pipe $j\ra i$, then $t_j\ra s_i$ is a pipe in $B$ and so $t_j = w(s_i) = v(i) = t_{u(i)}$ which implies $j = u(i)$.
Hence $u(i)\ra i$ is a pipe in $\varphi_v^w(B)$ for each $1\leq i\leq m$, so $\varphi_v^w(B)\in \BPD(u)$.

\begin{lemma}
\label{lem:remov_min}
    Let $w\in S_n$, let $v\sseq w$, and set $u:= \perm(v)$.
    Then $\varphi_v^w(B) \in \mBPD(u)$.
\end{lemma}

\begin{proof}
    The only part left to show is that $\varphi_v^w(B)$ is minimal.
    This is not difficult to see for if $j\ra i$ were a removable pipe in $\varphi_v^w(B)$, then $t_j\ra s_i$ is a pipe in $B$ which we claim is removable.
    Indeed, first note that $B_{s_i,t_j} = \varphi_v^w(B)_{i,j} = \relbow$.
    Now consider row $s_i$ in $B$.
    We know that $B_{s_i,t_k} = \varphi_v^w(B)_{i,k} \neq \relbow$ for all $1\leq k\leq m$, $k\neq j$.
    Furthermore, $B_{s_i,y_k}$ cannot be a $\relbow$ for any $1\leq k\leq r$ since otherwise pipe $y_k\ra x_k$ would not be removable in $B$.
    Hence $B_{s_i,t_j}$ is the only $\relbow$ in row $s_i$ and a similar argument shows the same for column $t_j$.
    On the other hand, $t_j\ra s_i$ cannot be removable since $y_1\ra x_1,\dots,y_r\ra x_r$ were assumed to be the only removable pipes in $B$, so we get a contradiction.
    Therefore $\varphi_v^w(B)$ must be minimal.
\end{proof}

% \begin{corollary}
% \label{cor:det}
%     $\det\Phi_v^w(B)$ agrees with $\det B$ up to a sign.
% \end{corollary}

\begin{proposition}
\label{prop:BPD_bij}
    Let $w\in S_n$, let $v\sseq w$, and set $u:= \perm(v)$.
    Then $\varphi_v^w: \BPD(w;v) \ra \mBPD(u)$ is a bijection.
\end{proposition}

\begin{proof}
    Lemma \ref{lem:remov_min} tells us that $\varphi_v^w$ is a well-defined map from $\BPD(w;v)$ to $\mBPD(u)$.
    We shall show that $\varphi_v^w$ is a bijection by constructing an inverse map.
    
    For $B\in \mBPD(u)$, perform the following:
    \begin{enumerate}
        \item[($1'$)] Expand out the columns in $B$ by moving them to positions $t_1<\dots<t_m$ and then adding $\blank$ and $\horizontal$ tiles as appropriate in the gap positions between columns (i.e., in columns $y_1,\dots,y_r$) to preserve the linkage of pipes.
        \item[($2'$)] Expand out the rows in $B$ by moving them to positions $s_1<\dots<s_m$ and then adding $\blank$ and $\vertical$ tiles as appropriate in the gap positions between rows (i.e., in rows $x_1,\dots,x_r$) to preserve the linkage of pipes.
        \item[($3'$)] Add undrooped pipes $y_1\ra x_1,\dots,y_r\ra x_r$ to $B$.
        We are able to do this since each row $x_k$ only contains $\blank$ and $\vertical$ tiles while each column $y_k$ only contains $\blank$ and $\horizontal$ tiles ($1\leq k\leq r$).
    \end{enumerate}
    Let $\psi_v^w(B)$ denote the diagram obtained by this process.
    % See Figure ref{?} for an example.
    It is clear that $\psi_v^w(B)$ is a bumpless pipe dream of size $n$.
    As before, $B_{i,j} = \psi_v^w(B)_{s_i,t_j}$ for all $1\leq i,j\leq m$.
    Hence each pipe $j\ra i$ in $B$ gets sent to the corresponding pipe $t_j\ra s_i$ in $\psi_v^w(B)$.
    Since we add the additional pipes $y_1\ra x_1,\dots, y_r\ra x_r$, $\psi_v^w(B)$ in fact has permutation $w$.
    
    We need to show that $\psi_v^w(B) \in \BPD(w;v)$; that is, the removable pipes of $\psi_v^w(B)$ are precisely $y_1\ra x_1,\dots,y_r\ra x_r$.
    It is clear from construction that these pipes are indeed removable since for each $1\leq k\leq r$, the elbow $B_{x_k,y_k} = \relbow$ is the only $\relbow$ that gets added to row $x_k$ and column $y_k$.
    
    Now, consider the pipe $t_j\ra s_i$ in $\psi_v^w(B)$ for some $1\leq i,j\leq m$.
    Then $j\ra i$ is a pipe in $B$ which is not removable since $B$ is minimal.
    It follows that $B_{i,k}=\relbow$ for some $1\leq k\leq m$, $k\neq j$.
    % or $B_{k,j}=\relbow$ for some $1\leq k\leq m$, $k\neq i$.
    %If the former, then 
    Then $\psi_v^w(B)_{s_i,t_k}=B_{i,k}=\relbow$ which implies that $s_i\ra t_j$ is not removable in $\psi_v^w(B)$.
    % A similar story happens with the second case.
    We conclude that $t_j\ra s_i$ is not removable in $\psi_v^w(B)$, so in fact $\psi_v^w(B)\in \BPD(w;v)$.
    
    It is clear that the procedure in $\psi_v^w$ will undo the procedure in $\varphi_v^w$ and vice versa, so $\psi_v^w = (\varphi_v^w)^{-1}$ as desired.
\end{proof}

We now examine how $\varphi_v^w$ acts on reduced bumpless pipe dreams.
The first statement we will show is that $\varphi_v^w$ preserves the property of being reduced i.e., $\varphi_v^w(B)$ is reduced whenever $B\in \bpd(w;v)$ is reduced.

\begin{proposition}
\label{prop:bpd_include}
    Let $w\in S_n$, let $v\sseq w$, and set $u:= \perm(v)$.
    Then $\varphi_v^w(\bpd(w;v)) \sseq \mbpd(u)$.
\end{proposition}

\begin{proof}
    Consider two pipes $j\ra i$ and $j'\ra i'$ in $\varphi_v^w(B)$.
    It is easy to see that $j\ra i$ and $j'\ra i'$ cross in $\varphi_v^w(B)$ at some $\varphi_v^w(B)_{a,b} = \cross$ if and only if $B_{s_a,t_b} = \varphi_v^w(B)_{a,b} = \cross$ is a cross shared by the pipes $t_j\ra s_i$ and $t_{j'}\ra s_{i'}$ in $B$.
    Hence the number of times $j\ra i$ and $j'\ra i'$ cross in $\varphi_v^w(B)$ is equal to the number of times $t_j\ra s_i$ and $t_{j'}\ra s_{i'}$ cross in $B$,
    but this latter quantity is at most one since $B$ is reduced.
    We conclude that $\varphi_v^w(B)$ is reduced.
\end{proof}

The previous proposition gives us the inequality $\abs{\bpd(w;v)} \leq \abs{\mbpd(u)}$
since $\varphi_v^w$ is injective.
In a moment, we will give a criterion on $w$ for when this is an equality.
For now, we have the following.

\begin{corollary}
\label{cor:schub_ub}
    There is an upper bound \[
        \nu_w = \abs{\bpd(w)}
        \leq \sum_{u\leq w} \abs{\mbpd(u)}\cdot p_u(w).
    \]
\end{corollary}

\begin{proof}
    We have \[
        \abs{\bpd(w)}
        = \sum_{v\sseq w} \abs{\bpd(w;v)}
        \leq \sum_{v\sseq w} \abs{\mbpd(\perm(v))}
        = \sum_{u\leq w} \abs{\mbpd(u)} \cdot p_u(w).
    \]
\end{proof}

For $B\in \BPD(w;v)$, observe that at no point in applying $\varphi_v^w$ to $B$ do we alter the number of $\jelbow$ tiles in $B$.
It follows that $\abs{(\varphi_v^w(B))^{-1}(\jelbow)} = \abs{B^{-1}(\jelbow)}$.
When $B\in \bpd(w;v)$ is reduced, we can say more since $\wt^{(\beta)}(B) = (1+\beta)^{|B^{-1}(\smalljelbow)|}$ only depends on $\abs{B^{-1}(\jelbow)}$ (since $\abs{B^{-1}(\blank)} = \ell(w)$).
Hence we get the following result.

\begin{proposition}
\label{prop:bpd_weight}
    Let $w\in S_n$, let $v\sseq w$, and set $u:= \perm(v)$.
    Then $\wt^{(\beta)}(B) = \wt^{(\beta)}(\varphi_v^w(B))$ for all $B\in \bpd(w;v)$.
\end{proposition}

\begin{proof}
    Let $B\in \bpd(w;v)$ be reduced.
    Then $\varphi_v^w \in \mbpd(u)$ is reduced by Proposition \ref{prop:bpd_include}, so \[
        \wt^{(\beta)}(B) =
        (1+\beta)^{|B^{-1}(\smalljelbow)|} =
        (1+\beta)^{|(\varphi_v^w(B))^{-1}(\smalljelbow)|} =
        \wt^{(\beta)}(\varphi_v^w(B)).
    \]
\end{proof}

\subsection{1243-avoiding and Vexillary Permutations}
\label{sect:1243_2143}

We will now investigate properties of permutations that avoid the patterns $1243$ and $2143$.

\begin{proposition}
\label{prop:BPD_nonred_1243_2143}
    Let $B\in \BPD(n)$ be a bumpless pipe dream that is not reduced.
    \begin{enumerate}
        \item[(a)] The permutation of $B$, $w_B$, contains a $1243$ pattern or a $2143$ pattern.
        In particular, if $B$ has two pipes that cross a nonzero even number of times, then $w_B$ contains a $1243$ pattern, and if $B$ has two pipes that cross an odd number of times, then $w_B$ contains a $2143$ pattern.
        \item[(b)] The type of $B$, $\partial(B)$, contains a $2143$ pattern.
    \end{enumerate}
\end{proposition}

\begin{proof}
    Since $B$ is not reduced, it has two pipes $y\ra x$ and $y'\ra x'$ ($x< x'$) that cross at least twice.
    Let $(i_1,j_1)$ and $(i_2,j_2)$ be the locations of the final two $\cross$\,s shared by $y\ra x$ and $y'\ra x'$ that we encounter when traversing $y\ra x$ from south to east.
    Note that the $\cross$ at $(i_1,j_1)$ must lie strictly to the southwest of the $\cross$ at $(i_2,j_2)$, whereas both of these $\cross$\,s lie southeast of both $(x,y)$ and $(x',y')$.
    In particular $i_1>i_2\geq x'$ and $j_2>j_1\geq \max\{y,y'\}$.
    In order to exit in row $x<x'$, $y\ra x$ must pass underneath $y'\ra x'$ horizontally through $(i_1,j_1)$ and then emerge back above vertically through $(i_2,j_2)$, at which point $y\ra x$ will stay above $y'\ra x'$ for the remainder of its path.
    This means that $y\ra x$ must engage in a $\jelbow$ at some point $(i_1,k)$ to the east of $(i_1,j_1)$ (so $k>j_1$) in order to turn from east to north and eventually hit $(i_2,j_2)$.
    
    Of the $\jelbow$\,s in $B$ that lie weakly southeast of $(i_1,k)$, let $(a,b)$ be the position of one chosen to be maximally southeast (so $a\geq i_1$, $b\geq k$ and $B_{a',b'}\neq \jelbow$ for all $a'\geq a$, $b'\geq b$).
    Observe that the pipe entering in column $b$ must have its first $\relbow$ at some position $(c,b)$ strictly south of $(a,b)$ (so $c>a$), at which point it cannot engage in a $\jelbow$ by the maximality of $(a,b)$.
    Similarly, the pipe exiting in row $a$ must have its final $\relbow$ at some position $(a,d)$ strictly to the east of $(a,b)$ (so $d>b$), below which no $\jelbow$ can occur since $(a,b)$ was chosen maximally.
    Hence $c\ra b$ and $d\ra a$ are each (undrooped) pipes in $B$.
    See Figure \ref{fig:1243_diagram} for a diagram.

\begin{figure}[h]
    \centering
    %
    % DIAGRAM 1
    %
    \begin{tikzpicture}[scale=0.5]
    \draw[step=10.0,black,opacity=0.5,thin,xshift=0.5cm,yshift=0.5cm] (0,0) grid (10,10);
    % crosses (i_1,j_1), (i_2,j_2)
    \node at (11,5) {$i_1$};
    \node at (11,8) {$i_2$};
    \node at (4,0) {$j_1$};
    \node at (5,0) {$k$};
    \node at (6,0) {$j_2$};
    \draw[gray]
    (4,5) circle (0.5)
    (6,8) circle (0.5);
    \filldraw[gray,opacity=0.2]
    (4.5,5.5) rectangle (10.5,0.5);
    % pipe y -> x
    \draw[color=\pipecolor,very thick]
    (1,0.5) to (1,3.75) to[bend left] (1.25,4) to (2.75,4) to[bend right] (3,4.25) to (3,4.75) to[bend left] (3.25,5) to (4.75,5) to[bend right] (5,5.25) to (5,6.75) to[bend left] (5.25,7) to (5.75,7) to[bend right] (6,7.25) to (6,9.75) to[bend left] (6.25,10) to (10.5,10);
    \node at (1,0) {$y$};
    \node at (11,10) {$x$};
    % pipe y' -> x'
    \draw[color=\pipecolor,very thick]
    (3,0.5) to (3,2.75) to[bend left] (3.25,3) to (3.75,3) to[bend right] (4,3.25) to (4,7.75) to[bend left] (4.25,8) to (6.75,8) to[bend right] (7,8.25) to (7,8.75) to[bend left] (7.25,9) to (10.5,9);
    \node at (3,0) {$y'$};
    \node at (11,9) {$x'$};
    % SE maximal elbows
    \draw[color=UFblue,thick]
    (4.5,2) to (4.75,2) to[bend right] (5,2.25) to (5,2.5)  % #1
    (6.5,3) to (6.75,3) to[bend right] (7,3.25) to (7,3.5)  % #2 (a,b)
    (8.5,4) to (8.75,4) to[bend right] (9,4.25) to (9,4.5); % #3
    \draw[color=UFblue,very thick,dotted]
    (3.5,2) to (4.5,2) (5,2.5) to (5,3.5)  % #1
    (5.5,3) to (6.5,3) (7,3.5) to (7,4.5)  % #2 (a,b)
    (7.5,4) to (8.5,4) (9,4.5) to (9,5.5); % #3
    % pipe d -> a
    \draw[color=\pipecolor,very thick]
    (9,0.5) to (9,2.75) to[bend left] (9.25,3) to (10.5,3);
    \node at (9,0) {$d$};
    \node at (11,3) {$a$};
    % pipe b -> c
    \draw[color=\pipecolor,very thick]
    (7,0.5) to (7,1.75) to[bend left] (7.25,2) to (10.5,2);
    \node at (7,0) {$b$};
    \node at (11,2) {$c$};
    %c\foreach \i\j in {1/10,3/9,7/2,9/3}{\node at (\i,\j) {$\bullet$};}
    \end{tikzpicture}
    \qquad
    %
    % DIAGRAM 2
    %
    \begin{tikzpicture}[scale=0.5]
    \draw[step=10.0,black,opacity=0.5,thin,xshift=0.5cm,yshift=0.5cm] (0,0) grid (10,10);
    % crosses (i_1,j_1), (i_2,j_2)
    \node at (11,5) {$i_1$};
    \node at (11,8) {$i_2$};
    \node at (4,0) {$j_1$};
    \node at (5,0) {$k$};
    \node at (6,0) {$j_2$};
    \draw[gray]
    (4,5) circle (0.5)
    (6,8) circle (0.5);
    \filldraw[gray,opacity=0.2]
    (4.5,5.5) rectangle (10.5,0.5);
    % pipe y -> x
    \draw[color=\pipecolor,very thick]
    (2,0.5) to (2,3.75) to[bend left] (2.25,4) to (2.75,4) to[bend right] (3,4.25) to (3,4.75) to[bend left] (3.25,5) to (4.75,5) to[bend right] (5,5.25) to (5,6.75) to[bend left] (5.25,7) to (5.75,7) to[bend right] (6,7.25) to (6,9.75) to[bend left] (6.25,10) to (10.5,10);
    \node at (2,0) {$y$};
    \node at (11,10) {$x$};
    % pipe y' -> x'
    \draw[color=\pipecolor,very thick]
    (1,0.5) to (1,2.75) to[bend left] (1.25,3) to (3.75,3) to[bend right] (4,3.25) to (4,7.75) to[bend left] (4.25,8) to (6.75,8) to[bend right] (7,8.25) to (7,8.75) to[bend left] (7.25,9) to (10.5,9);
    \node at (1,0) {$y'$};
    \node at (11,9) {$x'$};
    % SE maximal elbows
    \draw[color=UFblue,thick]
    (4.5,2) to (4.75,2) to[bend right] (5,2.25) to (5,2.5)  % #1
    (6.5,3) to (6.75,3) to[bend right] (7,3.25) to (7,3.5)  % #2 (a,b)
    (8.5,4) to (8.75,4) to[bend right] (9,4.25) to (9,4.5); % #3
    \draw[color=UFblue,very thick,dotted]
    (3.5,2) to (4.5,2) (5,2.5) to (5,3.5)  % #1
    (5.5,3) to (6.5,3) (7,3.5) to (7,4.5)  % #2 (a,b)
    (7.5,4) to (8.5,4) (9,4.5) to (9,5.5); % #3
    % pipe d -> a
    \draw[color=\pipecolor,very thick]
    (9,0.5) to (9,2.75) to[bend left] (9.25,3) to (10.5,3);
    \node at (9,0) {$d$};
    \node at (11,3) {$a$};
    % pipe b -> c
    \draw[color=\pipecolor,very thick]
    (7,0.5) to (7,1.75) to[bend left] (7.25,2) to (10.5,2);
    \node at (7,0) {$b$};
    \node at (11,2) {$c$};
    %\foreach \i\j in {1/9,2/10,7/2,9/3}{\node at (\i,\j) {$\bullet$};}
    \end{tikzpicture}
    \caption{The two configurations of the pipes $y\ra x$ and $y'\ra x'$ in $B$ as in the proof of Proposition \ref{prop:BPD_nonred_1243_2143} depending on if they cross an even or an odd number of times.
    The final two crosses shared by these pipes are encircled.
    In each diagram, the shaded region consists of the positions weakly southeast of the j-elbow at $(i_1,k)$ and are populated with candidate maximal southeast j-elbows (in blue).
    In the even configuration, we get the pattern $1243$; in the odd, the pattern $2143$.}
    \label{fig:1243_diagram}
\end{figure}
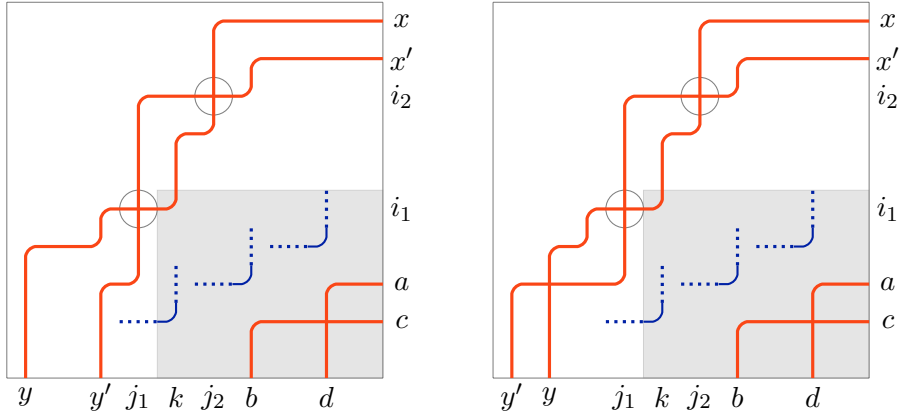

    We now have the sequences
    $x < x'\leq i_2 < i_1 \leq a < c$ and
    $\max\{y, y'\} \leq j_1 < k \leq b < d$.
    Since $w_B(s_i) = t_j$, $w_B(x_k) = y_k$, $w_B(a) = d$, and $w_B(c) = b$, $yy'db$ is a subword of $w$ that is either an occurrence of the pattern $1243$ (if $y<y'$) or the pattern $2143$ (if $y>y'$).
    The first situation occurs precisely when $y\ra x$ and $y'\ra x'$ cross an \emph{even} number of times, in which case $y\ra x$ enters and exists on the outside of $y'\ra x'$.
    When $y\ra x$ and $y'\ra x'$ cross an \emph{odd} number of times, $y\ra x$ still exits above $y'\ra x'$ but instead enters underneath $y'\ra x'$, so we get the second case.
    This completes (a).
    
    Observe that the undrooped pipes $b\ra c$ and $d\ra a$ lie to the southeast of the $\jelbow$ at $(a,b)$.
    The southeast maximality of $(a,b)$ tells us that there can be no $\jelbow$\,s underneath the hooks of either $b\ra c$ and $d\ra a$.
    In particular, $b\ra c$ and $d\ra a$ each cannot engage in double crosses with other pipes in $B$, so $b\ra c$ and $d\ra a$ are also pipes in $B_K$.
    
    Since the tile at $(i_1,j_1)$ is a $\cross$ in $B$, it is either a $\cross$ or a $\bump$ in $B_K$.
    Suppose that $t\ra s$ and $t'\ra s'$ ($s< s'$) are the pipes in $B_K$ meeting at $(i_1,j_1)$.
    These pipes must cross in $B_K$; this is obvious if $(B_K)_{i_1,j_1} = \cross$.
    If instead $(B_K)_{i_1,j_1}=\bump$, then these pipes must have crossed somewhere previously in $B_K$ as dictated by the construction of $B_K$.
    Since pipes in $B_K$ can cross at most once, this must be the case for $t\ra s$ and $t'\ra s'$. 
    Hence $t'\ra s'$ enters $B_K$ from the south on the outside of $t\ra s$ (since it exits to the east on the inside of $t\ra s$),  so $t'< t$.
    
    Since the pipes $t\ra s$ and $t'\ra s'$ visit the tile $(i_1,j_1)$, we must have, in particular, that $s'\leq i_1$ and $t\leq j_1$.
    Then $s< s'\leq i_1\leq a< c$ and $t'<t\leq i_1< k\leq b< d$.
    To finish, we note that $s'\neq a$ since $t'\ra s'$ and $d\ra a$ are pipes in $B_K$ with $t'\neq d$.
    Hence $tt'bd$ is an occurrence of the pattern $2143$ in $\partial(B)$ which gives (b).
\end{proof}

% Implicit in the above proof is the following corollary.

% \begin{corollary}
% \label{cor:BPD_nonred_1243_2143}
%     Let $w\in S_n$ and let $B\in \BPD(w)$.
%     \begin{enumerate}
%         \item[(a)] If $B$ has two pipes that cross a nonzero even amount of times, then $w$ contains a $1243$ pattern.
%         \item[(b)] If $B$ has two pipes that cross an odd amount of times, then $w$ contains a $2143$ pattern.
%     \end{enumerate}
% \end{corollary}

\begin{corollary}
\label{cor:BPD_1243_2143}
    Let $w\in S_n$.
    \begin{enumerate}
        \item[(a)] If $w$ is vexillary, then $\BPD_K(w) = \bpd(w)$.
        \item[(b)] If $w$ is vexillary and $1243$-avoiding, then 
        $\BPD(w) = \BPD_K(w) = \bpd(w)$.
    \end{enumerate}
\end{corollary}

\begin{proof}
    First suppose that $w$ is vexillary (i.e., $w$ is $2143$-avoiding).
    We already saw that in general $\bpd(w) \sseq \BPD(w)$ and $\bpd(w) \sseq \BPD_K(w)$.
    Proposition \ref{prop:BPD_nonred_1243_2143} shows that if there is a nonreduced $B\in \BPD_K(w)$, then $w$ must contain a $2143$ pattern.
    Since reduced bumpless pipe dreams have permutation $w$ if and only if they have type $w$,
    we conclude that $\BPD_K(w) = \bpd(w)$ which is (a).
    
    Now suppose additionally that $w$ avoids the pattern $1243$.
    If there is $B\in \BPD(w)$ nonreduced, then $w$ contains either $1243$ or $2143$ by Proposition \ref{prop:BPD_nonred_1243_2143}, so this cannot be the case.
    Hence $\BPD(w) = \bpd(w)$.
    Combine this result with (a) to get (b).
\end{proof}

We remark that Corollary \ref{cor:BPD_1243_2143}(a) was also proved by Weigandt \cite{Wei21}.

Now, assume the notation from the previous subsection.
In other words, for $w\in S_n$ a permutation, $v\sseq w$ a subword of size $m$, and $u:= \perm(v) \in S_m$, we denote the indices and entries of $v$ by $s_1<\cdots<s_m$ and $t_1<\cdots<t_m$, respectively, we assume each bumpless pipe dream $B\in \BPD(w;v)$ has removable pipes $y_1\ra x_1,\dots,y_r\ra x_r$ ($r:= n-m$, $y_1<\cdots<y_r$), and we take $\varphi_v^w: \BPD(w;v)\ra \mBPD(u)$ to be the bijection as in Proposition \ref{prop:BPD_bij}.

\begin{proposition}
\label{prop:bpd_bij_1243}
    Let $w\in S_n$, let $v\sseq w$, and set $u:= \perm(v)$.
    If $w$ is $1243$-avoiding, then
    $\varphi_v^w(\bpd(w;v)) = \mbpd(u)$.
    In this case, $\varphi_v^w$ restricts to a bijection $\bpd(w;v)\ra \mbpd(u)$.
\end{proposition}

\begin{proof}
    We will demonstrate that $\psi_v^w(\mbpd(u)) \sseq \bpd(w;v)$ (here $\psi_v^w$ is the inverse map to $\varphi_v^w$ we constructed in the proof of Proposition \ref{prop:BPD_bij}) in which case $\mbpd(u) = \varphi_v^w(\psi_v^w(\mbpd(u)) \sseq \varphi_v^w(\bpd(w;v))$.
    The result then follows by combining with Proposition \ref{prop:bpd_include}.
    
    Let $B\in \mbpd(u)$.
    Our goal is to show that $\psi_v^w(B)$ is reduced.
    To do so, we first note that there are two types of pipes in $\psi_v^w(B)$: pipes of the form $t_j\ra s_i$ where $1\leq i,j\leq m$, and pipes of the form $y_k\ra x_k$ for $1\leq k\leq r$.
    We will check case-wise that no two pipes of any combination of these types can cross twice in $\psi_v^w(B)$.
    
    Consider two pipes of first type, say $t_{j_1}\ra s_{i_1}$ and $t_{j_2}\ra s_{i_2}$.
    As in the proof of Proposition \ref{prop:bpd_include},
    it is clear that the number of times these pipes cross in $\psi_v^w(B)$ is the same as the number of times the pipes $j_1\ra i_1$ and $j_2\ra i_2$ cross in $B$.
    Since $B$ is reduced, it follows that $t_{j_1}\ra s_{i_1}$ and $t_{j_2}\ra s_{i_2}$ can cross at most once.
    
    We now consider two pipes of the second type, say $y_{k_1}\ra x_{k_1}$ and $y_{k_2}\ra x_{k_2}$.
    The observation that pipes of this type are undrooped in $\psi_v^w(B)$ is enough to check that no configuration of $y_{k_1}\ra x_{k_1}$ and $y_{k_2}\ra x_{k_2}$ will allow for multiple crosses.
    Thus we are done with this case as well.
    
    It remains to check that a pipe $t_j\ra s_i$ and a pipe $y_k\ra x_k$ cannot cross twice in $\psi_v^w(B)$.
    Suppose otherwise.
    Since $y_k\ra x_k$ is undrooped, it can only cross another pipe at most twice.
    Hence $t_j\ra s_i$ and $y_k\ra x_k$ cross exactly twice, so by Proposition \ref{prop:BPD_nonred_1243_2143}(a) $w$ must contain a $1243$ pattern.
    This contradicts our hypothesis on $w$, so we get the final case we needed in our analysis.
    Therefore $\psi_v^w(B)$ is reduced, so the proposition follows.
\end{proof}

% The consequence of Proposition \ref{prop:bpd_bij_1243} is that $\varphi_v^w$ restricts to a bijection between $\bpd(w;v)$ and $\mbpd(u)$ when $w$ is $1243$-avoiding.

\begin{proposition}
\label{prop:schub_1243}
    Let $w\in S_n$.
    If $w$ is $1243$-avoiding,
    then \[
        \nu_w = \abs{\bpd(w)}
        = \sum_{u\leq w} \abs{\mbpd(u)}\cdot p_u(w).
    \]
\end{proposition}

\begin{proof}
    Using Proposition \ref{prop:bpd_bij_1243}, we can replace the inequality with an equality in the proof of Corollary \ref{cor:schub_ub}.
\end{proof}

Observe that $|\mbpd(\emptyset)| = 1$ whereas $|\mbpd(w)| = \nu_w - \sum_{u<w} \abs{\mbpd(u)}\cdot p_u(w)$ by Proposition \ref{prop:schub_1243}.
These are the same seed and recurrence satisfied by the coefficients $c_w$ as defined in the introduction, so $c_w = |\mbpd(w)|$ for every $1243$-avoiding permutation $w$.
Theorem \ref{thm:gao_1243} immediately follows.

\begin{proposition}
\label{prop:bpd_weight_1243}
    Let $w\in S_n$, let $v\sseq w$, and set $u := \perm(v)$.
    If $w$ is $1243$-avoiding, then $\wt^{(\beta)}(\bpd(w;v)) = \wt^{(\beta)}(\mbpd(u))$.
\end{proposition}

\begin{proof}
    A simple calculation shows \begin{align*}
        \wt^{(\beta)}(\bpd(w;v))
        &= \sum_{B\in \bpd(w;v)} \wt^{(\beta)}(B) \\
        &= \sum_{B\in \bpd(w;v)} \wt^{(\beta)}(\varphi_v^w(B))
        \tag{Proposition \ref{prop:bpd_weight}} \\
        &= \sum_{B'\in \mbpd(u)} \wt^{(\beta)}(B')
        \tag{Proposition \ref{prop:bpd_bij_1243}} \\
        &= \wt^{(\beta)}(\mbpd(u)).
    \end{align*}
\end{proof}

\begin{proposition}
\label{prop:groth_1243_2143}
    Let $w\in S_n$.
    If $w$ is vexillary and $1243$-avoiding,
    then \[
        \nu^{(\beta)}_w = \wt^{(\beta)}(\BPD_K(w))
        = \sum_{u\leq w} \wt^{(\beta)}(\mbpd(u))\cdot p_u(w).
    \]
\end{proposition}

\begin{proof}
    We have
    \begin{align*}
        \wt^{(\beta)}(\BPD_K(w))
        &= \wt^{(\beta)}(\bpd(w))
        \tag{Corollary \ref{cor:BPD_1243_2143}(a)} \\
        &= \sum_{v\sseq w} \wt^{(\beta)}(\bpd(w;v)) \\
        &= \sum_{v\sseq w} \wt^{(\beta)}(\mbpd(\perm(v)))
        \tag{Proposition \ref{prop:bpd_weight_1243}} \\
        &= \sum_{u\leq w} \wt^{(\beta)}(\mbpd(u))\cdot p_u(w).
    \end{align*}
\end{proof}

As above, it is clear that $c^{(\beta)}_w = \wt^{(\beta)}(\mbpd(w))$ for all vexillary $1243$-avoiding permutations $w$.
This gives Theorem \ref{thm:groth_gao_1243_2143}.
Note that in this case, we also have $\mBPD(w) = \mbpd(w)$ by Corollary \ref{cor:BPD_1243_2143}(b), so in fact $c^{(\beta)}_w = \wt^{(\beta)}(\mBPD(w))$ as well.
\section{Closing Remarks}

\subsection{Future Directions}

The establishment of Theorem \ref{thm:gao_1243} cements $1243$-avoiding permutations as the second known (nontrivial) pattern-avoidance class of permutations for which Conjecture \ref{conj:gao} holds, the first being the avoidance class of $1432$ and $1423$ \cite{MT21}.
What is novel about this development is that it is done bijectively through Proposition \ref{prop:bpd_bij_1243} which, for $w$ a $1243$-avoiding permutation, essentially says we can embed a distinct copy of $\mbpd(u)$ into $\bpd(w)$ for each occurrence of the pattern $u$ in $w$.
In contrast, Conjecture \ref{conj:gao} was shown \cite{MT21} to hold for permutations $w$ avoiding $1432$ and $1423$ indirectly by verifying the identity \[
    \sum_{v\sseq w} (-1)^{\abs{w}-\abs{v}}\nu_{\perm(v)} \geq 0.
\]
With the observation $\nu_w = \sum_{v\sseq w} c_{\perm(v)}$, we see by inclusion-exclusion that $c_w$ is equal to the alternating sum on the left-hand side.

We may wonder if the methods used in this paper can be used to show Conjecture \ref{conj:gao} for other pattern-avoidance classes of permutations (including for permutations avoiding both $1432$ and $1423$).
% Of course, the ultimate goal would be to prove Conjecture \ref{conj:gao} in its full generality, but for now we may ponder if there are other pattern-avoidance classes of permutations for which the method outlined in Section 3 can apply.
One prerequisite is narrowing our definition of minimal bumpless pipe dreams.
To see why this is necessary, we examine the permutation $1243$ itself.
Recall that $\bpd(1243) = \BPD_K(1243)$ is given in Figure \ref{fig:red_bpd} (in orange).
Write $B_1$, $B_2$, and $B_3$ for the three diagrams in $\bpd(1243)$ labeled from left-to-right.
$1243$ has three subwords $v$ for which $c_{\perm(v)}\neq 0$: $143$ and $243$ (occurrences of $132$), and the empty word $\emptyset$.
Since $c_\emptyset = c_{132} = 1$, we would like to have $\bpd(1243,\emptyset)$, $\bpd(1243,143)$, and $\bpd(1243,243)$ to each contain a single BPD.
This is the case for $\bpd(1243,\emptyset) = \{B_1\}$ and $\bpd(1243,243) = \{B_2\}$, but not for $\bpd(1243,143)$ which is empty.
Instead $\mbpd(1243) = \{B_3\}$ has one too many bumpless pipe dreams (since $c_{1243} = \nu_w - c_\emptyset - c_{132}\cdot 2 = 3 - 1 - 1\cdot 2 = 0$).
Thus if we want our approach to work for $1243$, $B_3$ should no longer be considered minimal and instead be reclassified as part of $\bpd(1243,143)$.
This entails generalizing the notion of removable pipes (i.e., pipe $2\ra 2$ in $B_3$ should be removable) and adapting the removal map $\varphi_v^w$ to account for these extra pipes.

On the other hand, the above approach will not work for $1243$ in the more general setting of Conjecture \ref{conj:groth_gao} without more drastic changes to our method.
To understand why, we first calculate that $c^{(\beta)}_{1243} = \beta^2 + \beta$.
If our approach is to work, this should be the same as the $\beta$-weight of some collection of ``minimal'' bumpless pipe dreams of type $1243$.
Here we are doomed since no subset of $\BPD_K(1243)$ has this particular weighting.

To this end, it may prove advantageous to forgo the bumpless pipe dream model in favor of its older sibling: the \emph{pipe dream} model.
Pipe dreams are a family of wiring diagrams that provide an alternative description of Schubert and Grothendieck polynomials.
They were first introduced as \emph{RC-graphs} by Bergerson and Billey \cite{BB93} in reference to their identification with \emph{C}ompatible sequences of \emph{R}educed words, which are known to give a sum formula for Schubert polynomials \cite{BJS93}.
Pipe dreams are fundamentally different objects from bumpless pipe dreams and are amenable to different techniques of analysis.
In fact, the original $132$ and $132$-$1432$ bounds on $\nu_w$ of Weigandt \cite{Wei18} and Gao \cite{Gao21} were proved using the pipe dream model.
There may be a class of permutations for which the methods from this paper could be applied to pipe dreams, developing notions of removable pipes and minimality to this alternate setting.

\subsection{Layered Permutations}

Of separate interest is a question regarding the maximum values attained by $\nu_w$ and $c_w$ among $w\in S_n$ for a fixed size $n$.
It has been conjectured by Merzon and Smirnov \cite{MS16} that the permutations $w\in S_n$ maximizing $\nu_w$ are always \emph{layered} i.e., of the form $w = n_1 \dots 1\,
    n_2\dots (n_1+1)
    \ .\ .\ .\ 
    n_k\dots (n_{k-1}+1)
$ for positive integers $1\leq n_1< n_2< \cdots< n_k=n$, which has been confirmed for $n\leq 10$.
Gao \cite{Gao21} found that the coefficients $c_w$ are also maximized for layered $w\in S_n$, $n\leq 8$.

Continuing this trend, we have checked that for $n\leq 9$, the permutations $w\in S_n$ maximizing the $\beta=1$ specializations $\nu^{(1)}_w$ and $c^{(1)}_w$ are all layered.
Interestingly, the permutations maximizing $\nu^{(1)}_w$ are precisely those that maximize $c^{(1)}_w$.
This was not the case with the sets of permutations maximizing $\nu_w$ and $c_w$, which differ for $n=5$ \cite{Gao21}.
We have provided the following table (Figure \ref{fig:max}) which records our findings.

\begin{figure}[h]
    \centering
    \begin{tabular}{|c|l|l|l|}
        \hline
        $n$ & max $\nu^{(1)}_w$ & max $c^{(1)}_w$ & $w\in S_n$ \\
        \hline
        $0$ & $1$ & $1$ & $\emptyset$ \\
        $1$ & $1$ & $0$ & $1$ \\
        $2$ & $1$ & $0$ & $12$ and $21$ \\
        $3$ & $3$ & $2$ & $132$ \\
        $4$ & $11$ & $4$ & $1432$ \\
        $5$ & $71$ & $44$ & $12543$ and $21543$ \\
        $6$ & $1101$ & $828$ & $132654$ \\
        $7$ & $38259$ & $32160$ & $1327654$ \\
        $8$ & $1711251$ & $1501128$ & $13287654$ \\
        $9$ & $190013835$ & $177205856$ & $143298765$ \\
        \hline
    \end{tabular}
    \caption{The maximum values of $\nu^{(1)}_w$ and $c^{(1)}_w$ among all $w\in S_n$, as well as the $w\in S_n$ achieving these maxima.}
    \label{fig:max}
\end{figure}

The permutations listed above are in general not the same as the permutations maximizing $\nu_w$.
In particular, they differ for $n=5$, $n=6$, and $n=9$.
Likewise, they differ from the permutations maximizing $c_w$ for $n=6$ and $n=9$.\footnote{
    We checked for $n=9$ that $109294$ is the maximum value for $c_w$ which is achieved by $w = 132987654$.
    This is the same permutation maximizing $\nu_w$.
}
Except for $n=2$ and $n=5$, there is a unique permutation maximizing $\nu^{(1)}_w$ and $c^{(1)}_w$.
In the cases where we do not have uniqueness, the polynomials $\nu^{(\beta)}_w$ and $c^{(\beta)}_w$ agree for all permutations attaining the maximum.

\subsection{Formulae for Skew Sums}

We finish with a simple observation which hastens the computation of $\nu^{(\beta)}$ and $c^{(\beta)}$ for skew sums of permutations.
Recall that the \emph{skew sum} of two permutations $u\in S_m$ and $v\in S_n$ is the permutation $u\ominus v\in S_{m+n}$ given by $(u(1)+n)\dots (u(m)+n)\, v(1)\dots v(n)$.
When is $(i,j)$ an inversion of $u\ominus v$?
If $i,j\in[m]$, then this is precisely when $(i,j)$ is an inversion of $u$,
and likewise if $i,j\in [m+1,m+n]$, then this happens when $(i-m,j-m)$ is an inversion of $v$.
Finally, if $i\in[m]$ and $j\in[m+1,m+n]$, then $(i,j)$ is always an inversion of $u\ominus v$.
We thus see that $\ell(u\ominus v) = m\cdot n + \ell(u) + \ell(v)$.

We claim that the diagrams in $\BPD_K(u\ominus v)$ are precisely those of the form
\begin{equation}
    \label{eq:skew_sum}
    B\ =\ \begin{pmatrix}
    \blank_{m\times n} & B_u \\
    B_v & \cross_{n\times m}
    \end{pmatrix}
    % \begin{tabular}{|c|c|}
    % \hline
    % $\blank_{m\times n}$ & $B_u$ \\
    % \hline
    % $B_v$ & $\cross_{n\times m}$ \\
    % \hline
    % \end{tabular}
\end{equation}
where $B_u\in \BPD_K(u)$, $B_v\in \BPD_K(v)$, $\blank_{m\times n}$ is an $m\times n$ grid of $\blank$ tiles, and $\cross_{n\times m}$ is an $n\times m$ grid of $\cross$ tiles.
To see this, let $B\in \BPD_K(u\ominus v)$.
There are two types of pipes in $B_K$: pipes of the form $y\ra x$ with $x\in [m]$ (which correspond to $u$) and pipes of the form $y'\ra x'$ with $x'\in [m+1,m+n] := \{m+1,\dots,m+n\}$ (which correspond to $v$).
Each pipe $y\ra x$ of the first type is restricted to the region $[m+n]\times [n+1,m+n]$ since $y = u\ominus v(x) = u(x)+n \geq n+1$.
Similarly, each pipe $y'\ra x'$ of the second type is restricted to $[m+1,m+n]\times[m+n]$ since $x\geq m+1$.
In particular, no pipe in $B_K$ can enter the region $[m]\times [n]$, so it must be all $\blank$ tiles.
Hence $B$ must also be $\blank$\,s in the region $[m]\times [n]$.

Furthermore, every pair of pipes $y\ra x$ and $y'\ra x'$ in $B_K$ with $x\in [m]$, $x'\in [m+1,m+n]$ must cross at least once since $(x,x')$ is an inversion in $u\ominus v$.
This can only happen in the region $[m+1,m+n]\times [n+1,m+n]$ in $B_K$.
On the other hand, there are $m\cdot n$ such pairs, so every cell in this region must be a $\cross$, each being the only cross belonging to a unique pair $y\ra x$ and $y'\ra x'$.
Thus $B$ is also $\cross$\,s in the region $[m+1,m+n] \times [n+1,m+n]$.
From here it is easy to see that the other two regions in $B$ must be bumpless pipe dreams of type $u$ and type $v$, respectively, so $B$ has the desired form.

Conversely, suppose that $B$ is a diagram of the above form.
Observe that each $\cross$ in the region $\cross_{n\times m}$ in $B$ belongs to a unique pair of pipes $y\ra x$ and $y'\ra x'$ where $x'\in [m+1,m+n]$, $x\in [m]$.
This fact remains unchanged no matter how many $\cross$\,s in $B_u$ or $B_v$ we change into $\bump$\,s (though the specific pipes passing through each $\cross$ in $\cross_{n\times m}$ may change).
Hence when passing to $B_K$, the region $\cross_{n\times m}$ remains unchanged.
From this it is easy to see that \[
    B_K\ =\ \begin{pmatrix}
    \blank_{m\times n} & (B_u)_K \\
    (B_v)_K & \cross_{n\times m}
    \end{pmatrix},
\]
so $B$ has type $u\ominus v$.
This completes the claim.

For $B\in \BPD_K(u\ominus v)$, write $B = B_u\ominus B_v$ to denote the decomposition of $B$ as in Equation \ref{eq:skew_sum}.
% with corresponding sub-diagrams $B_u\in \BPD_K(u)$ and $B_v\in \BPD_K(v)$ as above.
% Write $B = B_u\ominus B_v$ to denote this relationship.
Clearly $\abs{B^{-1}(\blank)} = m\cdot n + \abs{B_u^{-1}(\blank)} + \abs{B_v^{-1}(\blank)}$
and $\abs{B^{-1}(\jelbow)} =  \abs{B_u^{-1}(\jelbow)} + \abs{B_v^{-1}(\jelbow)}$.
Coupled with our previous computation that $\ell(u\ominus v) = m\cdot n + \ell(u) + \ell(v)$, it easily follows that $\wt^{(\beta)}(B) = \wt^{(\beta)}(B_u) \wt^{(\beta)}(B_v)$.

We can now compute
\begin{align*}
    \nu^{(\beta)}_{u\ominus v}
    &= \sum_{B\in \BPD_K(u\ominus v)} \wt^{(\beta)}(B) \\
    &= \sum_{B_u\in \BPD_K(u)} \sum_{B_v\in \BPD_K(v)} \wt^{(\beta)}(B_u\ominus B_v) \\
    &= \sum_{B_u\in \BPD_K(u)} \sum_{B_v\in \BPD_K(v)} \wt^{(\beta)}(B_u) \wt^{(\beta)}(B_v) \\
    &= \Big(\sum_{B_u\in \BPD_K(u)} \wt^{(\beta)}(B_u)\Big)
    \Big(\sum_{B_v\in \BPD_K(u)} \wt^{(\beta)}(B_v)\Big)
    = \nu^{(\beta)}_u\nu^{(\beta)}_v.
\end{align*}
This gives the first multiplicative formula in the following proposition.

\begin{proposition}
\label{prop:skew_formulae}
    Let $u\in S_m$ and $v\in S_n$.
    Then $\nu^{(\beta)}_{u\ominus v} = \nu^{(\beta)}_u \nu^{(\beta)}_v$ and 
    $c^{(\beta)}_{u\ominus v} = c^{(\beta)}_u c^{(\beta)}_v$.
\end{proposition}

\begin{proof}
It remains to show that $c^{(\beta)}_{u\ominus v} = c^{(\beta)}_u c^{(\beta)}_v$.
The subwords of $u\ominus v$ are themselves of the form $u'\ominus_n v'$ for subwords $u'\sseq u$ and $v'\sseq v$,
where by $u'\ominus_n v'$ we mean the word $(u'(1)+n)\dots (u'(|u'|)+n)\ v'(1)\cdots v'(|v'|)$.
It is clear that $\perm(u'\ominus_n v') = \perm(u')\ominus \perm(v')$.
% Note that the permutation patterns of $u\ominus v$
% are themselves just skew sums $u'\ominus v'$ of permutations patterns $u'\leq u$ and $v'\leq v$.
% The number of times that $u'\ominus v'$ appears as a pattern in $u\ominus v$ in this way is given by $p_{u'}(u) p_{v'}(v)$.
By applying inclusion-exclusion, we then get \begin{align*}
    c^{(\beta)}_{u\ominus v}
    &= \sum_{u'\ominus_n v'\, \sseq\, u\ominus v}
    (-1)^{\abs{u\ominus v}-|u'\ominus v'|}
    \nu^{(\beta)}_{\perm(u'\ominus v')} \\
    &= \sum_{u'\sseq u} \sum_{v'\sseq v}
    (-1)^{(\abs{u}+\abs{v}) - (|u'|+|v'|)}
    \nu^{(\beta)}_{\perm(u')\ominus \perm(v')} \\
    &= \sum_{u'\sseq u} \sum_{v'\sseq v}
    (-1)^{\abs{u}-|u'|}
    (-1)^{\abs{v}-|v'|}
    \nu^{(\beta)}_{\perm(u')}
    \nu^{(\beta)}_{\perm(v')} \\
    &= \Big(\sum_{u'\sseq u}
    (-1)^{\abs{u}-|u'|}
    \nu^{(\beta)}_{u'}\Big)
    \Big(\sum_{v'\sseq v}
    (-1)^{\abs{v}-|v'|}
    \nu^{(\beta)}_{v'}\Big) \\
    &= c^{(\beta)}_u c^{(\beta)}_v.
\end{align*}
\end{proof}

\section*{Acknowledgement}

The author is grateful to Zachary Hamaker, who has served as their mentor throughout each step of this project.
Special thanks is given to Adam Gregory, whose conversations with the author helped direct research toward this topic and make early progress.
This work %was performed while the author was an undergraduate at the University of Florida and
was partially supported by NSF Grant DMS-2054423.

\end{document}